\documentclass[a4paper,12pt]{amsart}
\usepackage{amsmath}
\usepackage{amssymb}
\usepackage{amsfonts}

\def \endproof {\quad \hfill  \rule{2mm}{2mm} \par\medskip}
\renewcommand{\theequation}{\thesection.\arabic{equation}}
\DeclareMathSymbol{\subsetneqq}{\mathbin}{AMSb}{36}

\def \endproof {\quad \hfill  \rule{2mm}{2mm} \par\medskip}
\DeclareMathSymbol{\subsetneqq}{\mathbin}{AMSb}{36}

\newcommand{\R}{\mathbb{R}}

\newcommand{\dint}{\displaystyle\int}
\newcommand{\dsum}{\displaystyle\sum}

\newtheorem{thm}{{\bf Theorem}}[section]
\newtheorem{lem}{{\bf Lemma}}[section]
\newtheorem{prop}{{\bf Proposition}}[section]

\theoremstyle{remark}
\newtheorem{rem}{\bf Remark}[section]

\theoremstyle{definition}
\newtheorem{defi}{\bf Definition}[section]

\author{Mohamed Majdoub}
\address{Facult\'e des Sciences de Tunis, D\'epartement
de Math\'ematiques\\ Campus universitaire 1060, Tunis, Tunisia.}
\email{\it mohamed.majdoub@fst.rnu.tn}
\thanks{M.M is partially supported by the {\sf Laboratory of PDE and applications}
of Faculty of Sciences, Tunis, Tunisia}
\author{Marius Paicu}
\address{Universit\'e de Paris-Sud, D\'epartement de Math\'ematiques\\
B\^{a}timent 425, Orsay, France }

\email{\it marius.paicu@math.u-psud.fr}
\thanks{}

\subjclass{..., ..., ....} 

\title{Uniform Local Existence for Inhomogeneous Rotating Fluid Equations}
\date{\today}

\begin{document}
\begin{abstract}
We investigate the equations of anisotropic incompressible viscous
fluids in $\R^3$, rotating around an inhomogeneous vector
  $B( t, x_1, x_2)$. We prove the global existence of
 strong solutions in suitable anisotropic Sobolev spaces for small initial data, as well as uniform
 local existence result with respect to the Rossby number in the
 same functional spaces under the additional assumption that
 $B=B(t,x_1)$ or $B=B(t,x_2)$. We also obtain the propagation of
 the isotropic Sobolev regularity using a new refined product law.
\end{abstract}

\maketitle

\tableofcontents \vspace{ -1\baselineskip}

\renewcommand{\theequation}{\thesection.\arabic{equation}}

\newpage


\section{Introduction}

The aim of this work is to investigate the well-posedness theory of
inhomogeneous rotating fluids in the framework of Sobolev spaces. We
consider the equations governing incompressible, viscous fluids in
$\R^3$, rotating around an inhomogeneous vector
$B(t,x_h)=\big(b_1(t,x_h),b_2(t,x_h),b_3(t,x_h)\big)$ where $x_h$
stands for the horizontal variables (that is $x=(x_h,x_3)$). This is
a generalization of the usual rotating fluid model, where
$B=e_3:=(0,0,1)$ . More precisely, we are interested in the
following system
\begin{eqnarray*}({IRF^\varepsilon})
\begin{cases}
\partial_t u^\varepsilon+ (u^\varepsilon.\nabla) u^\varepsilon-\nu_h\Delta_h
u^\varepsilon-\nu_v\partial_{x_3}^2u^\varepsilon+\dfrac{1}{\varepsilon}\, \Big(u^\varepsilon\times
B\Big)+\nabla p^\varepsilon=0\quad\mbox{in}\quad \R^3\\
\mbox{div}\,u^\varepsilon=0\quad\mbox{in}\quad \R^3\\
u^\varepsilon(0,x)=u_0(x)
\end{cases}
\end{eqnarray*}
where $u$ is the velocity field and $p$ is the pressure. The
constants $\nu_h>0$, $\nu_v\geq 0$ and $\varepsilon>0$ represent
respectively the horizontal viscosity, the vertical viscosity and
the Rossby number. We have written
$\Delta_h=\partial_{1}^2+\partial_{2}^2$  and
$(u.\nabla)u=\sum_{j=1}^3 u_j\partial_j u$ where
$u=(u_1,u_2,u_3)$ and $\partial_j$ represents the partial derivative
in the direction $x_j$. We will also write
$\nabla_h=(\partial_1,\partial_2)$ and $\mbox{div}_h\, u=\partial_1
u_1+\partial_2 u_2$. Throughout this article, $B$ will be a smooth
function with bounded derivatives. Additional assumptions on $B$
will be made later.

Before going any further, let us recall some well-known facts about
the constant case $B=e_3$.  There is a large literature dealing with
these equations in the constant case.  The first results for systems
with large skew-symmetric perturbations are contained in \cite{KM},
\cite{JMR1}, \cite{JMR2}.  In \cite{BMN}, A. Babin, A. Mahalov and B.
Nicolaenko studied the incompressible rotating Euler and Navier-Stokes
equations on the torus.  Using the method introduced by S. Schochet
 (see \cite{S1} and \cite{S2}), I. Gallagher studied in \cite{G} and
\cite{G1} this problem in its abstract hyperbolic form.  In the case
of the incompressible rotating Navier-Stokes equations on the torus, it
is shown (see \cite{BMN} and \cite{Gr}) that the solutions exist
globally in time and converge to a solution of some diffusion
equation.  The case of rotating fluids evolving between two horizontal
plates with Dirichlet boundary conditions was studied by
E. Grenier and N. Masmoudi \cite{GM} (see also the paper of N.
Masmoudi \cite{Mas1} and the recent work of J.-Y. Chemin, B.
Desjardin, I. Gallagher and E. Grenier \cite{CDGG1}).  Motivated by
this case, J.-Y. Chemin {\it{and al.}} studied the incompressible
fluids with anisotropic viscosity in the whole space in \cite{CDGG}
and \cite{CDGG1}.  They obtained local existence of the solution in
the anisotropic Sobolev spaces $H^{0,1/2+\eta}(\R^3)$ and the global
existence for data which are small compared to the horizontal
 viscosity.  They also proved the global existence of the solution for
anisotropic rotating fluids.  The key of their proof is an anisotropic
version of Strichartz estimates.  In \cite{paicu1},
global existence of the solution for rotating fluids with vanishing
vertical viscosity was shown in the periodic case.  Finally, we also refer
to the related works \cite{BMN1}, \cite{bresch}, \cite{Cheverry} and
\cite{GSR1}.

In both constant and non constant cases, if $\nu_v>0$, it is easy to
construct global weak solutions.  This is simply due to the fact that
the singular perturbation is a linear skew-symmetric operator.  The
behavior of these solutions (as $\varepsilon\rightarrow 0$) has been
studied by I. Gallagher and L. Saint-Raymond in a recent paper
\cite{GR}.  It is proved in \cite{GR} that the weak solutions
$u^\varepsilon$ converge to the solution of a heat equation when the
vector field $B$ is non stationary.  The proofs of \cite{GR} are based
on weak compactness arguments.  For a detailed analysis of the
rotating geophysical fluids we refer to \cite{GSR}.  \\

In this paper, we would like to investigate the question of
existence of strong solutions for (${IRF^\varepsilon}$). Our aim is
to prove the existence and uniqueness of a solution on a uniform
time interval, or the global existence and uniqueness for small
initial data. To do so we need an energy estimate in Sobolev spaces as in \cite{FK}.
Unfortunately one sees quickly that this is not an easy matter just
by taking horizontal derivatives of the equation. But since $B$ does
not depend on the third variable, all vertical derivatives are
allowed in the energy estimate. Only horizontal derivatives create
extra terms. So the idea is to start by taking the initial data in
an anisotropic-type Sobolev space. Such spaces have been introduced
by Iftimie in \cite{Iftimie} and are very useful for anisotropic
problems (see also \cite{CDGG},
  \cite{Iftimie1, Iftimie2}, \cite{paicu, paicu1}). Using the skew-symmetry of
the Coriolis operator in the anisotropic Sobolev space $H^{0,s}$
(see Definition \ref{Hss'} below) and some product laws we shall
prove the global existence and uniqueness for small data in
$H^{0,s}$ for $s>1/2$.

The next step consists in proving uniform local existence and
uniqueness for large data in the same anisotropic Sobolev space. For
technical reasons, we are not able to prove such result in the
general case but only under the assumption that the rotation vector
depends on one space variable. Let us remark that this assumption on
$B$ is consistent with some models of geophysical flows (see
\cite{GSR}).

Once these steps are achieved, we return to the propagation of the
isotropic Sobolev regularity. This is a delicate problem du to the
lack of vertical viscosity. The major difficulty is to estimate the
term $\left\langle (u.\nabla) u,u\right\rangle_{H^s}$ in a good way.
In the third section we obtain a new refined product law (see Lemma
\ref{s}) which enables
  us to propagate the $H^s$ regularity as
claimed in Theorems \ref{thNS}, \ref{global} and \ref{th3}.


\subsection{Statement of the results}

Our first result is the global existence of solutions of
$(IRF^\varepsilon)$ in
  suitable anisotropic Sobolev spaces when the
  initial data are small enough.

   \begin{thm}
\label{th1} \, \, Assume that $B=B(t,x_1,x_2)$. Let $s>1/2$ be a real number and  $u_0$
be a divergence-free vector field in $H^{0,s}(\R^3)$. Assume that
$\|u_0\|_{H^{0,s}}\leq c\,\nu_h$ where $c$ is small enough. Then,
there exists a unique global solution $u^\varepsilon$ of
($IRF^\varepsilon$) such that
$$u^\varepsilon\in C_b(\R_+; H^{0,s})\hskip0.3cm\text{and}\hskip0.3cm \nabla_h u^\varepsilon\in L^2(\R_+; H^{0,s}).$$
Moreover, for any $t\geq 0$,
$$\|u^\varepsilon(t)\|_{H^{0,s}}^2+2\nu_h
\dint_0^t\|\nabla_h u^\varepsilon(\tau)\|_{H^{0,s}}^2\,d\tau\leq
\|u_0\|_{H^
 {0,s}}^2.
$$
\end{thm}

\smallskip

 The proof of this theorem involves no real difficulties
if one takes into account the earlier results of \cite{paicu}. It is
based on the following simple fact
$$
\langle u\times B(t,x_h), u\rangle_{H^{0,s}}=0.
$$

The next theorem gives the local existence, uniformly in
$\varepsilon$, in the space $H^{0,s}$ when the rotation vector only
depends on one space variable, say $x_1$.
\smallskip

\begin{thm}
\label{th2} Assume that $B=B(t,x_1)$. Let $s>1/2$ be a real number
and $u_0\in H^{0,s}(\R^3)$ be a divergence-free vector field. Then,
there exists a positive time $T$ such that, for all $\varepsilon>0$,
the system ($IRF^\varepsilon$) has a unique solution $u^\varepsilon$
such that
$$u^\varepsilon\in {\mathcal C}([0,T], H^{0,s})\hskip0.3cm\text{and}
\hskip0.3cm\nabla_h u^\varepsilon\in L^2([0,T], H^{0,s}).$$
\end{thm}

The assumption on $B$ made in this theorem seems to be physical.
  Indeed, in the literature, $B$ usually  depends only on
 one
variable, for example $B=(1+\beta x_2)e_3$, where $x_2$ is the latitude
variable and $\beta$ is a parameter. Let us explain briefly the mathematical reasons for the assumption
made on $B$.  Generally speaking, we can prove local existence of the solution
as long as we can prove that the solution of the linear part of the
equation stays uniformly small in a short time. In the general case
($B=B(t,x_1,x_2)$), we can prove only that the solution of the
linear corresponding equation is bounded, but we cannot prove that
 the solution is small for short time. Contrary to this situation, if $B=B(t,x_1)$ or $B=B(t,x_2)$ we can prove that  the solution of the linear problem is small in uniformly small time and so we can prove the uniformly local existence result.

Theorem \ref{th2} is not completely satisfactory, since the initial
data belongs to an anisotropic Sobolev space which is adapted to the
equation. A natural question to ask is then the following: is it
possible to propagate the isotropic Sobolev regularity?\\

Let us first consider the particular case $B\equiv 0$, that is, the case  of the anisotropic Navier-Stokes equations.

\begin{eqnarray*}(NS_h)
\begin{cases}
\partial_t u+ u.\nabla u-\nu_h\Delta_h
u+\nabla p=0\\
\mbox{div}\,u=0\\
u(0,x)=u_0(x)
\end{cases}
\end{eqnarray*}

This system was studied in \cite{CDGG}, where the local existence for
arbitrary initial data and global existence for small initial data
in $H^{0,s}(\R^3)$ for $s>1/2$ and uniqueness for $s>3/2$ are
proved. In \cite{Iftimie2}, the author filled the gap between
existence and uniqueness and proved that uniqueness holds when
existence does; that is for $s>1/2$. In the third section we will
prove the following theorem which can be seen as a propagation of
the isotropic Sobolev regularity in time.

\begin{thm}
\label{thNS} Let $s>1/2$ be a real number, and let $u_0\in
H^s(\R^3)$ be a divergence-free vector field. There exists a
positive time $T$ and a unique solution $
 u$ of ($NS_h$) defined on
$[0,T]\times\R^3$  such that
$$u\in {\mathcal C}([0,T]; H^s) \quad\mbox{and}\quad \nabla_h u\in
L^2([0,T]; H^s).$$ \noindent If the maximal time $T^*$ of existence
is finite, then
 $$\lim_{t \rightarrow
T^*}\int_0^t\|\nabla_h
u(\tau)\|_{L^\infty_v(L^2_h)}^2\big(1+\|u(\tau)\|^2_{L^\infty_v(L^2_h)}\big)d\tau=+\infty.$$
Furthermore, there exists a positive constant $c$ such that if
$\|u_0\|_{H^s}\leq c\, \nu_h$, then the solution is global in time.
\end{thm}

The fact that the lack of vertical viscosity prevents us from
gaining vertical regularity is the main difficulty in the proof of
Theorem \ref{thNS}. The main tool to overcome this serious
difficulty is a new refined product law stated in  Section
 3 (see Lemma \ref{s}) and proved in Section 6.

From this theorem we derive the following two results. The first one
concerns the propagation of the isotropic Sobolev regularity for
($IRF^\varepsilon$). The second one concerns the uniform local existence
in the isotropic Sobolev space.
\begin{thm}
\label{global} Let $s>1/2$ be a real number. There exists $c>0$
small enough such that, for any  divergence-free vector field
$u_0\in H^s(\R^3)$ with $\|u_0\|_{H^s}\leq c\nu_h$, and for any
$\varepsilon>0$, the system ($IRF^\varepsilon$) has a unique global
solution $u^\varepsilon$ such that
$$
u^\varepsilon\in{\mathcal C}(\R_+; H^s(\R^3))\cap
L^\infty_{loc}(\R_+;H^s(\R^3))\quad\mbox{and}\quad\nabla_h
u^\varepsilon\in L^2_{loc}(\R_+; H^s(\R^3)).
$$
\end{thm}

\smallskip

\begin{thm}
\label{th3} Assume that $B=B (t,x_1)$. Let $s>1/2$ be a real number,
and let $u_0\in H^s(\R^3)$ be a divergence-free vector field. Then
,there exists a positive time $T$ such that, for all
$\varepsilon>0$, the system ($IRF^\varepsilon$) has a unique
solution $u^\varepsilon$ satisfying
 $$u^\varepsilon\in{\mathcal
C}([0,T]; H^s)\quad\mbox{and}\quad \nabla_h u^\varepsilon\in
L^2([0,T]; H^s).$$
\end{thm}

\begin{rem}
 We point out that all the results
  stated above still holds in the periodic
case. The proofs are similar except for some technical differences
(see \cite{paicu2}).
  \end{rem}

  \smallskip

  The structure of the paper is as follows: Section 2 contains some notations
  and technical results which will be used in the whole paper. The third section deals
   with the propagation of the isotropic Sobolev regularity for the anisotropic
   Navier-Stokes equations. Section four  is devoted to the
  proofs of global existence results for small initial data. In the fifth section, we give  the proofs of the uniform local existence
  theorems. Section six deals with the proof of the main product law stated in the third section (see Lemma \ref{s}). Finally,  the appendix is devoted to the proofs
  of some technical lemmas and product laws, which
  are more or less  contained in earlier papers.
  For the sake of completeness we give their proofs here.


\section{Notations and technical lemmas}

In view of the anisotropy of the problem, we shall have to use
functions spaces that take into account this anisotropy. The main
tools will be the energy estimate and anisotopic Sobolev spaces.
Such spaces have been introduced by Iftimie in \cite{Iftimie} and
used by several authors (see for instance \cite{paicu,
paicu1}-\cite{CDGG}-\cite{Iftimie1, Iftimie2}).

The definition of these spaces requires an anisotropic dyadic
decomposition in the Fourier spaces. Let us first recall the
isotropic Littlewood-Paley theory. We refer to \cite{chem1} for a
precise decomposition of the Fourier space.


\subsection{Isotropic Littlewood-Paley theory}

The main idea is to localize in frequencies. The interest of this
theory consists in the fact that it allows us to define, at least
formally, the product of two distributions as  para-products
 and remainder.\\

Consider two smooth, compactly supported functions $\varphi$ and
$\chi$, with support respectively in a fixed ring of $\R$ far from
the origin, and in a fixed ball containing the origin and such
that
   $$ \forall\,\, s\in
\R\backslash\{0\},\quad\sum_{j\in\mathbb{Z}} \varphi
(2^{-j}s)=1\quad\mbox{and}\quad  \forall\,\, s\in \R,\quad
\chi(s)+\sum_{j\in\mathbb{N}} \varphi (2^{-j}s)=1$$

 \noindent Let us note that there exists an integer $N_0$ such that if
 $|j-k|\geq N_0$, then
 $\mbox{supp}\,\varphi(2^{-j}\cdot)\cap\varphi(2^{-k}\cdot)=\emptyset$.\\

 \noindent Next, we define the following operators of localization in Fourier
space:

\begin{eqnarray*}
{\text{ if}}\,\, j\geq 0,\,\,{\mathcal F}(\Delta _j u)(\xi
)&=&\varphi (2^{-j}|\xi|) {\mathcal F}(u)(\xi)\\ {\mathcal F}(\Delta
_{-1}u)(\xi )&=&\chi (|\xi|) {\mathcal
  F}(u) (\xi )
\\\Delta_j u&=&0\quad\mbox{for}\quad j<-1
\end{eqnarray*}

 \noindent where
$${\mathcal
F}(u)(\xi)=\hat{u}(\xi)=\dint_{\mathbb{R}^3}\;e^{-i\;x.\xi}\;u(x)\;dx,\quad
\xi=(\xi_1,\xi_2,\xi_3):=(\xi_h,\xi_3)\in\mathbb{R}^3 $$ is the
Fourier transform of the function $u$.\\

We recall also the definition of the anisotropic Lebesgue spaces:
 \begin{defi}
  \label{Lebesgue}
We denote by $L^p_h(L^r_v)$ the space
  $L^p(\R_{x_1}\times\R_{x_2};L^r(\R_{x_3}))$ endowed with the norm
  $$\|f\|_{L^p_h(L^r_v)} =\left\|\|f(x_h,.)\|_{L^r_{x_3}}\right\|_{L^p_{x_1, x_2}}.$$
  In the same way, we define $L^r_v(L^p_h)$ to be the space
   $L^r(\R_{x_3};L^p(\R_{x_1}\times\R_{x_2}))$ endowed with the norm
   $$
   \|f\|_{L^r_v(L^p_h)}=\left\|\|f(.,x_3)\|_{L^p_{x_1, x_2}}\right\|_{L^r_{x_3}}.$$
  \end{defi}

\smallskip

 The operators $\Delta_j$ satisfy the following property:

\begin{lem}
\label{uniformBound} The operator $\Delta_j$ is uniformly bounded in
$L^p_vL^r_h$   for all $1\le
 q p\leq\infty$ and $1\leq r\leq\infty$.
  \end{lem}

  \begin{lem}
  There exists a positive constant $C$ such that, for every $a\in\mathcal S'(\R^3)$ with
  $\nabla a\in L^\infty_v(L^s_h)$ and $b\in L^p_v(L^t_h)$\quad\\
\label{commutator}
$$
\Big\|[\Delta_j; a]b\Big\|_{L^p_v L^r_h}\leq C 2^{-j} \|\nabla
a\|_{L^\infty_v L^s_h} \|b\|_{L^p_v L^t _h}
$$
where $1/r=1/s +1/t$, $1\leq p\leq\infty$ and $[\Delta_j;
a]b=\Delta_j(a b)-a\Delta_j b$.
\end{lem}
\begin{lem}
\label{anablaa} There exists a positive constant $C$ such that, for
every $a\in\mathcal S'(\R^3)$, for $j\geq 0$, $1\leq p\leq\infty$ and
$1\leq
r\leq\infty$, we have\quad\\
$$
\|\Delta_j a\|_{L^p_v L^r_h}\leq C2^{-j}\|\Delta_j\nabla a\|_{L^p_v
L^r_h}.
$$
\end{lem}

The proofs of Lemmas \ref{uniformBound}, \ref{commutator} and
\ref{anablaa} are essentially
  contained in \cite{paicu} (see also \cite{chem} for a more general setting).
  For the sake of completeness we give the proofs in the Appendix.

  \

We recall now Bony's decomposition  \cite{bony}. It is well known that the
dyadic decomposition is useful for
defining the product of two distributions. Formally, we can write,
for any two distributions $u$ and $v$
$$
u=\dsum_{q}\,\Delta_q\,u\quad; \quad
v=\dsum_{q}\,\Delta_q\,v \quad\mbox {and}\quad
uv=\dsum_{q, q'}\, \Delta_q\,u\,\Delta_{q'}\,v.
$$
We split this sum into three terms
$$uv=T_uv+T_vu+R(u,v),$$
where we have denoted
 \begin{eqnarray*}
  T_u\,v&=&\dsum_{q'\leq q-2}\,\Delta_{q'} u\,\Delta_q v=\dsum_{q}\,
  S_{q-1} u\,\Delta_q v\\\\
R(u,v)&=&\dsum_{q}\dsum_{j\in\{0,\pm 1\}}\,\Delta_q u\,\Delta_{q+j}v.
\end{eqnarray*}
In the first term, the
 frequencies of $u$ are smaller than those of $v$. In the
second term, the frequencies of $v$ are smaller than
those of $u$. In the third term, the frequencies
of $u$ and $v$ are comparable. The first two  sums are called the
para-products and the third sum is the remaining term.

\noindent We have the following quasi-orthogonality properties

\begin{eqnarray*}
\Delta_q\left(S_{q'-1}u\,\Delta_{q'} v\right)&=&0\quad{\mbox if} \quad|q-q'|\geq N_0\\\\
\Delta_q \left(S_{q'+1} u\,\Delta_{q'} v\right)&=&0\quad{\mbox if} \quad q'\leq q-N_0.
\end{eqnarray*}

\noindent Note that $u$ is in the isotropic Sobolev space
$H^{s}(\R^3)$ if and only if
$$
\Big(\displaystyle\sum_{q}2^{2qs} \|\Delta_q
u\|^2_{L^2}\Big)^{1/2}<\infty,
$$
with equivalence of the  norms.


\subsection{The anisotropic case}

 In the frame of anisotropic Lebesgue spaces, the H\"older and Young
inequalities reads.
\begin{prop}
\label{Hol.You} We have the following inequalities:
$$
  \|f\,g\|_{L^r_v(L^p_h)}\leq\, \|f\|_{L^{r'}_v(L^{p'}_h)}
\|g\|_{L^{r''}_v(L^{p''}_h)}
  $$
  where $1/r=1/r'+1/r''$ and $1/p=1/p'+1/p''$.
  $$
\|f\star
g\|_{L^r_v(L^p_h)}\leq\,\|f\|_{L^{r'}_v(L^{p'}_h)}\|g\|_{L^{r''}_v(L^{p''}_h)}
$$
where $1+1/r=1/r'+1/r''$ and $1+1/p=1/p'+1/
 p''$.
\end{prop}

We recall the definition of anisotropic Sobolev spaces and their properties.
  \begin{defi}
  \label{Hss'}
  Let $s$ and $s'$ be two real numbers. The anisotropic Sobolev space
$H^{s,s'}(\R^3)$
  is the space of tempered distributions $u$ such that $\hat{u}$
  belongs to $L^1_{loc}(\R^3)$ and
  $$
  \|u\|_{H^{s,s'}}^2:=\dint_{\R^3}\,\left(1+|\xi_h|^2\right)^s\, \left(1+|\xi_3|^2\right)^{s'}|\hat{u}(\xi)|^2\,
  d\xi<\infty.
  $$
   \end{defi}

  \smallskip

  It is obvious that the space $H^{s,s'}(\R^3)$ endowed with the norm $\|{\bf
.}\|_{H^{s,s'}}$ is a Hilbert space.\\

As in the isotropic case, we have the following interpolation
property (see \cite{Iftimie}).
  \begin{prop}
  \label{interpolation}
Let $s,t,s',t'\in \R$ and $0\leq\alpha\leq 1$. If $u$ is in
$H^{s,s'}(\R^3)\cap H^{t,t'}(\R^3)$, then  $u$ belongs to $H^{\alpha
s+(1-\alpha)t,\alpha s'+(1-\alpha)t'}(\R^3)$ and
$$
\|u\|_{H^{\alpha s+(1-\alpha)t,\alpha s'+(1-\alpha)t'}}\,\leq
\|u\|_{H^{s,s'}
 }^{\alpha} \|u\|_{H^{t,t'}}^{1-\alpha}.
$$
  \end{prop}
  The multiplicative properties of the anisotropic Sobolev spaces
have been studied by
  several authors \cite{CDGG, Iftimie, Iftimie1, ST}.
  The following result is proved in \cite{Iftimie1}.
  \begin{thm}
  Let $s,t<1$, $s+t>0$, and $s', t'<\frac{1}{2}$, $s'+t'>0$. There exists
a constant $C>0$ such that, for any $u\in H^{s,s'}(\R^3)$ and
  $v\in H^{t,t'}(\R^3)$, the product  $u v$ belongs to $H^{s+t-1, s'+t'-1/2}(\R^3)$ and

  $$
\|u v\|_{H^{s+t-1, s'+t'-1/2}}\,\leq
C\|u\|_{H^{s,s'}}\,\|v\|_{H^{t,t'}}.
$$
\end{thm}
This theorem is not sufficient for our purpose since the regularity
we need is greater than $1/2$ in the vertical direction. The
following theorem, proved in \cite{Iftimie2}, deals with this
difficulty.
\begin{thm}
Let $s,t<1$, $s+t>0$, and $s'>\frac{1}{2}$,  $t'\leq s'$,
$s'+t'>0$. If $u\in H^{s,s'}(\R^3)$ and
  $v\in H^{t,t'}(\R^3)$, then $u v\in H^{s+t-1, t'}(\R^3)$ and there exists a
constant $C>0$ such that, for any $u\in H^{s,s'}(\R^3)$,  $v\in H^{t,t'}(\R^3)$, we have
  $$
\|u v\|_{H^{s+t-1, t'}(\R^3)}\,\leq
C\|u\|_{H^{s,s'}(\R^3)}\,\|v\|_{H^{t,t'}(\R^3)}.
$$
\end{thm}
\noindent The proofs of the above theorems and those of some product
laws, which will be stated below,  essentially use the anisotropic
Littlewood-Paley theory. Let us briefly recall  this theory (for
more details, we refer to  \cite{chem, paicu, paicu1}).\\

The basic idea of this theory is to localize in vertical
frequencies. Let us introduce
the operators of localization in vertical frequencies in the following way:

\begin{eqnarray*}
{\text{ if}}\,\, j\geq 0,\,\,\Delta_j^{\mbox v}\,u&=&{\mathcal
F}^{-1}\left(\varphi(2^{-j}|\xi_3|){\mathcal F}u(\xi)\right)\\
\Delta _{-1}^{\mbox v}u&=&{\mathcal F}^{-1}\left(\chi (|\xi_3|) {\mathcal F}(u)
(\xi )\right)\\
\Delta_j^{\mbox v} u&=&0\quad\mbox{for}\quad j<-1.
\end{eqnarray*}
We define also the operator
$$S_q^{\mbox v}\,u=\dsum_{q'\leq q-1}\,\Delta_{q'}^
 {\mbox v}\,u.$$

The interest of this decomposition consists in the fact that any
vertical derivative of a function localized in vertical frequencies
of size $2^q$ acts like the multiplication by $2^q$. The following
lemma, which is often referred to as Bernstein Lemma, will be very
useful in what follows.
\begin{lem}
\label{Bernstein} Let $u$ be a function such that
$\text{supp}\,{\mathcal F}^v\,u\subset \R^2_h\times 2^q{\mathcal
C}$, where ${\mathcal C}$ is a dyadic ring. Let $p\geq 1$ and $r\geq
r'\geq 1$ be real numbers. Then, there exists a constant $C>0$ such
that
\begin{equation}
\label{Berns1}
2^{qk}C^{-k}\|u\|_{L^p_h(L^r_v)}\leq\|\partial_{x_3}^k
u\|_{L^p_h(L^r_v)}\leq 2^{qk}C^{k}\|u\|_{L^p_h(L^r_v)}
\end{equation}
\begin{equation}
\label{Berns2}
2^{qk}C^{-k}\|u\|_{L^r_v(L^p_h)}\leq\|\partial_{x_3}^k
u\|_{L^r_v(L^p_h)}\leq 2^{qk}C^{k}\|u\|_{L^r_v(L^p_h)}
\end{equation}
\begin{equation}
\label{Berns3} \|u\|_{L^p_h(L^r_v)}\leq C
2^{q(\frac{1}{r'}-\frac{1}{r})}\|u\|_{L
 ^p_h(L^{r'}_v)}
\end{equation}
\begin{equation}
\label{Berns4} \|u\|_{L^r_v(L^p_h)}\leq C
2^{q(\frac{1}{r'}-\frac{1}{r})}\|u\|_{L^{r'}_v(L^p_h)}.
\end{equation}
\end{lem}

\

We give now the
Bony's decomposition in the vertical variable (see \cite{bony},
  \cite{CDGG}, \cite{Iftimie}). Formally, we can write,
for two distributions $u$ and $v$
$$
u=\dsum_{q}\,\Delta_q^{\mbox v}\,u\quad; \quad
v=\dsum_{q}\,\Delta_q^{\mbox v}\,v \quad\mbox {and}\quad
uv=\dsum_{q, q'}\, \Delta_q^{\mbox v}\,u\,\Delta_{q'}^{\mbox v}\,v.
$$
We split this sum into three terms. In the first term, the
vertical  frequencies of $u$ are smaller than those of $v$ while in the second term we have that   the vertical  frequencies of $v$ are smaller than
those of $u$. The first two  sums are called the
para-products . The third term is called the remaining term and is such that the vertical frequencies
of $u$ and $v$ are comparable.  We shall denote:
  \begin{eqnarray*}
  T_u\,v&=&\dsum_{q'\leq q-2}\,\Delta_{q'}^{\mbox v}u\,\Delta_q^{\mbox
v}v=\dsum_{q}\,
  S_{q-1}^{\mbox v}u\,\Delta_q^{\mbox v}v\\\\
R(u,v)&=&\dsum_{q}\dsum_{j\in\{0,\pm 1\}}\,\Delta_q^{\mbox
v}u\,\Delta_{q+j}^{\mbox v}v.
\end{eqnarray*}

\noindent We have the following quasi-orthogonality properties

\begin{eqnarray*}
\Delta_q^{\mbox v}\left(S_{q'-1}^{\mbox v}u\,\Delta_{q'}^{\mbox
v}v\right)&=&0\quad{\mbox if} \quad|q-q'|\geq N_0\\\\
\Delta_q^{\mbox v}\left(S_{q'+1}^{\mbox v}u\,\Delta_{q'}^{\mbox
v}v\right)&=&0\quad{\mbox if} \quad q'\leq q-N_0.
\end{eqnarray*}

\noindent Note that $u$ is in the anisotropic Sobolev space
$H^{0,s}(\R^3)$ if and only if
$$
\Big(\displaystyle\sum_{q}2^{2qs} \|\Delta_q^{\mbox v}
u\|^2_{L^2}\Big)^{1/2}<\infty,
$$
with equivalence of the  norms.


\section{Propagation of regularity  for the anisotropic
Navier-Stokes system}

In this section we will prove Theorem \ref{thNS}. Let us first
remark that this theorem is only a result about the propagation of
the isotropic Sobolev regularity of the initial data in the case of
anisotropic Navier-Stokes equations. The existence of the solution
is known in the larger space $H^{0,s}$ with $s>1/2$ (\cite{CDGG})
and in the critical Besov space $\mathcal B^{0,1/2}$ (\cite{paicu}).
We shall use the energy estimate in $H^s$ to prove that the solution
constructed in \cite{CDGG} preserves the isotropic Sobolev
regularity. The main difficulty is the control of the term
$\left\langle (u.\nabla) u,u\right\rangle_{H^s(\R^3)}$ without using
vertical derivatives of $u$. The following lemma is the key
 argument in the proof.
\begin{lem}
\label{s} Let $s>1/2$ be a real number. There exists a positive
constant $C$ such that, for any divergence-free vector field $u\in
H^s(\R^3)$, we have
\begin{eqnarray*}
\Big|\left\langle u.\nabla
u,u\right\rangle_{H^s}\Big|&\leq&C\Big[\|u\|_{L^\infty_v(
L^2_h)}^{1/2} \|\nabla_h u\|_{L^\infty_v
(L^2_h)
 }^{1/2}\|u\|_{H^s}^{1/2}\|\nabla_h u\|_{H^s}^{3/2}+\\&&
\|\nabla_h u\|_{L^\infty_v (L^2_h)}\|u\|_{H^s}\|\nabla_h
u\|_{H^s}\Big].
\end{eqnarray*}
\end{lem}

\noindent The proof of this lemma will be given in Section 6. Let us
go back to the proof of Theorem \ref{thNS}.

\begin{proof}[Proof of Theorem \ref{thNS}]\quad\\
\noindent The energy estimate in $H^s(\R^3)$ yields
\begin{equation}
\nonumber
\frac{d}{dt} \|u(t)\|_{H^s}^2+2\nu_h
\|\nabla_h u(t)\|_{H^s}^2\leq 2\Big|\langle u.\nabla u,
u\rangle_{H^s}\Big|.
\end{equation}
Using  Lemma \ref{s},  we obtain

\begin{eqnarray*}
\Big|\langle u.\nabla u,
u\rangle_{H^s}\Big|&\leq&C\Big[\|u\|_{L^\infty_v(L^2_h)}^{1/2}
\|\nabla_h u\|_{L^\infty_v(L^2_h)}^{1/2}\|u\|_{H^s}^{1/2}\|\nabla_h
u\|_{H^s}^{3/2}+\\&& \|\nabla_h
u\|_{L^\infty_v(L^2_h)}\|u\|_{H^s}\|\nabla_h u\|_{H^s}\Big]\\&\leq&
\nu_h/2 \|\nabla_h
u\|_{H^s}^2+C\,\nu_h^{-3}\Big(\nu_h^2+\|u\|_{L^\infty_v(L^2_h)}^2\Big)
\|\nabla_h u\|_{L^\infty_v(L^2_h)}^2\|u\|_{H^
 s}^2.
\end{eqnarray*}
The Gronwall's Lemma gives
\begin{equation}
\nonumber
\|u(t)\|_{H^s}^2\leq
\|u_0\|_{H^s}^2\exp\Big(\int_0^t\left(C\nu_h^{-3}(\nu_h^2+\|u(\tau)\|_{L^\infty_v(L^2_h)}^2)
\|\nabla_h u(\tau)\|_{L^\infty_v(L^2_h)}^2\right)d\tau\Big).
\end{equation}
This implies that, if
$$\int_0^t\left(C\nu_h^{-3}(\nu_h^2+\|u(\tau)\|_{L^\infty_v(L^2_h)}^2)
\|\nabla_h u(\tau)\|_{L^\infty_v(L^2_h)}^2\right)d\tau<+\infty,$$
for all $0\leq t\leq T$, then $T^*\geq T$, where $T^*$ is the
lifespan of the solution in $H^s$.

\noindent In particular, if $\|u_0\|_{H^s}\leq c\nu_h$ where $c$ is
small enough, then the anisotropic norm  of $u_0$ is small enough
($\|u_0\|_{H^{0,s}}\leq c\nu_h$) and  by the results proved in
\cite{CDGG}, there exists a unique global in time solution in the
space $H^{0,s}$. In addition
$$\forall\,\,t\geq 0,\,\,\,\|u(t)\|^2_{H^{0,s}}+\nu_h\int_0^t\|\nabla_h
u(\tau)\|^2_{H^{0,s}}d\tau\leq \|u_0\|_{H^{0,s}}^2.$$

\noindent Using the Sobolev embedding $H^{0,s}\hookrightarrow
L^\infty_v(L^2_h)$ for $s>1/2$, we get
$$\forall\,\,t\geq 0,\,\,\,\|u(t)\|^2_{L^\infty_v(L^2_h)}+\nu_h\int_0^t\|\nabla_h
u(\tau)\|^2_{L^\infty_v(L^2_h)}d\tau\leq \|u_0\|^2_{H^{0,s}}.$$
Therefore, the solution exists globally in $H^s$ and it is small
compared with the horizontal viscosity. This ends the proof of
Theorem \ref{thNS}.\end{proof}


\section{Global existence for small initial data}

In this section, we shall prove the global existence of solutions of
$(IRF^\varepsilon)$ in $H^{0,s}(\R^3)$ and $H^s(\R^3)$, when the
initial data are small enough. It should be noted here that the 2D
Navier-Stokes equations are globally well-posed in $L^2$, and that
the 3D case can be dealt with by splitting the horizontal and
vertical variables and using Sobolev embedding
$H^{0,s}\hookrightarrow L^\infty_v(L^2_h)$.
\\
The main tool in the proof of Theorem \ref{
 th1} is the following
product law (see \cite{CDGG} and \cite{paicu} for the critical
case).

\begin{lem}
\label{u,v,0,s} Let $s>1/2$ be a real number. There exists a
positive constant $C_s$ such that, for any  $u, v\in H^{0,s}(\R^3)$,
$u$ being divergence-free, we have
\begin{eqnarray*}
\Big|\langle u.\nabla v, v\rangle_{H^{0,s}}\Big|\leq
C_s\big(\|u\|^{1/2}_{H^{0,s}}\|\nabla_h
u\|^{1/2}_{H^{0,s}}\|v\|_{H^{0,s}}^{1/2}\|\nabla_h v\|_{H^{0,s}}^{3/2}\\
+\|v\|_{H^{0,s}}\|\nabla_h v\|_{H^{0,s}}\|\nabla_h
u\|_{H^{0,s}}\big).
\end{eqnarray*}
In particular
$$
\Big|\langle u.\nabla u, u\rangle_{H^{0,s}}\Big|\leq
C_s\|u\|_{H^{0,s}}\|\nabla_h u\|_{H^{0,s}}^2.
$$
\end{lem}
\noindent Note that, in such inequalities, the following simple fact
is fundamental
$$
\mbox{div}\,u=0\Longrightarrow\partial_3 u_3=-{\text div}_h u_h.
$$
\begin{proof}[Proof of Theorem \ref{th1}]\quad\\
Without loss of generality, one can assume that $\nu_v=0$. For
simplicity, we omit the dependance on $\varepsilon$
 in
$u^\varepsilon$. Therefore, $u$ satisfies
\begin{equation}
\label{NSTA}
\partial_t u+ (u .\nabla) u-\nu_h\Delta_h
u+\dfrac{1}{\varepsilon}\, \left(u\times B(t,x_h)\right)=-\nabla p.
\end{equation}
Taking the scalar product of (\ref{NSTA}) with $u$ in $H^{0,s}$, and
using the fact that
$$
\left\langle u\times B(t,x_h), u\right\rangle_{H^{0,s}}=0,
$$
we infer
\begin{equation}
\label{En1} \frac{d}{dt}\|u(t)\|_{H^{0,s}}^2+2\nu_h \|\nabla_h
u(t)\|_{H^{0,s}}^2\leq 2\Big|\langle (u.\nabla) u,
u\rangle_{H^{0,s}}\Big|.
\end{equation}
By Lemma \ref{u,v,0,s}, (\ref{En1}) becomes
\begin{equation}
\label{En2} \frac{d}{dt}\|u(t)\|_{H^{0,s}}^2+2\nu_h \|\nabla_h
u(t)\|_{H^{0,s}}^2\leq 2C_s\|u\|_{H^{0,s}}\|\nabla_h
u\|_{H^{0,s}}^2.
\end{equation}
Define $T^*$ by
$$
T^*=\sup\left\{\, T>0\,/\,\,\forall\,\, 0\leq t\leq
T,\quad\|u(t)\|_{H^{0,s}}\leq 2c\nu_h <\frac{\nu_h}{2C_s}\,\right\}.
$$
Using (\ref{En2}), we get
$$
\forall\,\,\, 0\leq
t<T^*,\,\,\,\|u(t)\|_{H^{0,s}}^2+\nu_h\,\dint_0^t \|\nabla_h
u(\tau)\|_{H^{0,s}}^2\,d\tau\leq \|u_0\|_{H^{0,s}}^2\leq (c\nu_h)^2.
$$
Thus
$$
\forall\,\,\, 0\leq t<T^*,\,\,\,\|u(t)\|_{H^{0,s}}\leq c\nu_h<2
c\nu_h,
$$
and so
$$
T^*=+\infty.
$$
Moreover, we have the following energy estimate
$$
\forall\,\,\, t\geq 0,\quad\|u(t)\|_{H^{0,s}}^2+\nu_h\,\dint_0^t
\|\nabla_h u(\tau)\|_{H^{0,s}}^2\,d\tau\leq \|u_0\|_{H^{0,s}}^2.
$$
Theorem \ref{th1} is then completely proved.\end{proof}

We now come to the proof of Theorem \ref{global} about global
existence for small initial data in the isotropic Sobolev space
$H^s$. Here, the product law stated in Lemma \ref{s} play a crucial
role.

\begin{proof}[Proof of Theorem \ref{global}]\quad\\
\noindent The energy estimate in $H^s(\R^3)$ implies
\begin{equation}
\nonumber
\frac{d}{dt}\|u(t)\|_{H^s}^2+2\nu_h \|\nabla_h u(t)\|_{H^s}^2\leq
2\Big|\langle (u.\nabla) u,
u\rangle_{H^s}\Big|+\frac{2}{\varepsilon}\Big|\langle u\times B,
u\rangle_{H^s}\Big|.
\end{equation}
The assumption on
 $B$ yields
$$
\Big|\langle u\times B, u\rangle_{H^s}\Big|\leq C\, \|u\|_{H^s}^2.
$$
Using Lemma \ref{s} and the injection $H^{0,s}\hookrightarrow
L^\infty_v(L^2_h)$, we infer
\begin{eqnarray*}
\Big|\langle (u.\nabla )u,
u\rangle_{H^s}\Big|&\leq&C\Big[\|u\|_{H^{0,s}}^{1/2} \|\nabla_h
u\|_{H^{0,s}}^{1/2}\|u\|_{H^s}^{1/2}\|\nabla_h u\|_{H^s}^{3/2}+\\&&
\|\nabla_h u\|_{H^{0,s}}\|u\|_{H^s}\|\nabla_h u\|_{H^s}\Big]\\&\leq&
\nu_h/2 \|\nabla_h
u\|_{H^s}^2+C\,\nu_h^{-3}\Big(\nu_h^2+\|u\|_{H^{0,s}}^2\Big)
\|\nabla_h u\|_{H^{0,s}}^2\|u\|_{H^s}^2.
\end{eqnarray*}
The Gronwall's Lemma leads to
\begin{equation}
\nonumber
\|u(t)\|_{H^s}^2\leq
\|u_0\|_{H^s}^2\exp\Big(\int_0^t\left(C/\varepsilon+
C\nu_h^{-3}(\nu_h^2+\|u(\tau)\|_{H^{0,s}}^2) \|\nabla_h
u(\tau)\|_{H^{0,s}}^2\right)d\tau\Big).
\end{equation}
The fact that $u_0$ is small in $H^{0,s}$ together with
Theorem~\ref{th1} implies that $u$ is global in time. Moreover
$$
\|u(t)\|_{H^{0,s}}^2+\nu_h\int_0^t\,\|\nabla_h
u(\tau
 )\|_{H^{0,s}}^2d\tau\leq \|u_0\|_{H^{0,s}}^2, \text{\;\;\; for
any\;\;\;}t\geq 0.
$$
It follows that
$$
\|u(t)\|_{H^s}^2\leq
\|u_0\|_{H^s}^2\exp\Big(C\nu_h^{-3}(\nu_h^2+\|u_0\|_{H^{0,s}}^2)
\|u_0\|_{H^{0,s}}^2\Big)\,\exp\Big(C\frac{t}{\varepsilon}\Big).
$$
This completes the proof of  Theorem \ref{global}.\end{proof}


\section{Uniform local existence}

In this section we wish to investigate the uniform local existence
in both anisotropic and isotropic cases. We prove Theorem \ref{th2}
and Theorem \ref{th3} stated in the introduction. As remarked
before, we shall assume that $B$ only depends on one space variable
since we are unable to recover similar results in general case.
Before we come to the details of the proofs, let us first
explain rapidly the scheme and make some clarifying comments.\\

We split the initial data $u_0$ in a small part in $H^{0,s}$ and a
smooth one. First, we
  solve the linear system with smooth initial
data to obtain a global and bounded solution $v^{\varepsilon}$.
Then, we consider the perturbed system satisfied by
$w^{\varepsilon}:=u^{\varepsilon}-v^{\varepsilon}$ and with small
initial data. We have to prove that $w^{\varepsilon}$ exists on a
uniform time interval. The strategy here is to use energy estimate
in $H^{0,s}$ and a Gronwall's lemma in order to prove that
$w^{\varepsilon}$ remains small on a uniform time interval. The main
difficulty comes from the term
$$
\langle (v^{\varepsilon}.\nabla )v^{\varepsilon},
w^{\varepsilon}\rangle_{H^{0,s}},
$$
but using the fact that two frequency directions of
$v^{\varepsilon}$ are bounded, one can prove
$$
\Big|\langle (v^{\varepsilon}.\nabla )v^{\varepsilon},
w^{\varepsilon}\rangle_{H^{0,s}}\Big|\leq\,C\|v^\varepsilon\|_{L^2}^3\|\nabla_h
v^{\varepsilon}\|_{L^2}+\frac{\nu_h}{10}\|\nabla_h
w^{\varepsilon}\|_{H^{0,s}}^2.
$$
This inequality, together with the fact that
$$
\dint_0^t\,\
 |v^\varepsilon(\tau)\|_{L^2}^3\|\nabla_h
v^{\varepsilon}(\tau)\|_{L^2}\,d\tau\lesssim\,t^{1/2}\,\|u_0\|_{L^2}^4,
$$
allows us to carry out the proofs.\\

We now come to the details.

\begin{proof}[Proof of Theorem \ref{th2}]\quad\\
As mentioned before, the main argument consists in splitting the
initial data in a small part in $H^{0,s}$ and a regular one. To do
this, we set $S_N u=\mathcal F^{-1}(\chi(2^{-N}|\xi|)\mathcal F u)$
and we take $N$ sufficiently large (depending on $u_0$ and $\nu_h$)
such that  $\|(I-S_N)u_0\|_{H^{0,s}}\leq c\nu_h$, where $c>0$ is
small enough. Then,   $u_0=S_N u_0+(I-S_N)u_0$, with  $S_Nu_0\in
H^{\infty}(\R^3)$. We point out that throughout the rest of this
section, $C$ will denotes a positive constant which may depends on
$N$, $\nu_h$, $s$  but not on $\varepsilon$.\\

\noindent First, we consider the following linear system:
\begin{equation}
\label{LS}
\begin{cases}
\partial_t v^{\varepsilon}_N-\nu_h\Delta_h
v^{\varepsilon}_N+\frac{1}{\varepsilon}\left(v^{\varepsilon}_N\times
B(t,x_1)\right)=-\nabla p^{\varepsilon}_N\\
\mbox{div}\,v^{\varepsilon}_N=0\\
v^{\varepsilon}_N|_{t=0}=S_Nu_0.
\end{cases}
\end{equation}

\noindent Clearly (\ref{LS}) has a unique global solution
$$v^{\varepsilon}_N\in C(\R_+;
H^{0,s})\hskip0.5cm\text{with}\hskip0.5cm\nabla_h
v^{\varepsilon}_N\in L^2(\R_+; H^{0,s}).$$

\noindent Since $B$ only depends on $t$ and $x_1$, we obtain
$S_N^{x_2,x_3}(u\times B(t,x_1))=S_N^{x_2,x_3}u\times B(t,x_1)$,
where $S_N^{x_2,x_3}$ is the operator of localization in $\xi_2$ and
$\xi_3$ direction defined by $S_N^{x_2,x_3}u={\mathcal
F}^{-1}\left({\tilde \chi}(2^{-N}|(\xi_2,\xi_3)|) {\mathcal F}
u(\xi)\right)$, $\tilde\chi$ being a smooth compactly supported
function whose value is $1$ on the support of $\chi$. By uniqueness
and the fact that $S_N^{x_2,x_3}S_Nu_0=S_Nu_0$, we get
$v^\varepsilon_N=S_N^{x_2,x_3}v^\varepsilon_N$. This implies  the
important fact that $v^\varepsilon_N$ is  a regular function with
respect to the $x_2$ and  $x_3$ variables.

\noindent Next, we consider the perturbed system:
\begin{equation}
\label{w}
\begin{cases}
\partial_t w^{\varepsilon}_N+(w^{\varepsilon}_N+
v^{\varepsilon}_N)\nabla (w^{\varepsilon}_N+
v^{\varepsilon}_N)-\nu_h\Delta_h w^{\varepsilon}_N+
\frac{1}{\varepsilon}\left(w^{\varepsilon}_N\times
B(t,x_1)\right)=-\nabla
q^{\varepsilon}_N\\
\mbox{div}\, w^{\varepsilon}_N=0\\
w^{\varepsilon}_N|_{t=0}=(I-S_N)u_0.
\end{cases}
\end{equation}
We shall now be interested in proving that $w^{\varepsilon}_N$ exists on a uniform time interval.
The energy estimate in  $H^{0,s}$ yields
\begin{eqnarray*}
\frac{1}{2}\frac{d}{dt}\|w^{\varepsilon}_N(t)\|^2_{H^{0,s}}+\nu_h
\|\nabla_hw^{\varepsilon}_N(t)\|^2_{H^{0,s}}\leq |\langle
w^{\varepsilon}_N\nabla
w^{\varepsilon}_N,w^{\varepsilon}_N\rangle_{H^{0,s}}|\\
+|\langle v^{\varepsilon}_N\nabla
w^{\varepsilon}_N,w^{\varepsilon}_N\rangle_{H^{0,s}}|+ |\langle
w^{\varepsilon}_N\nabla
v^{\varepsilon}_N,w^{\varepsilon
 }_N\rangle_{H^{0,s}}|\\
  +|\langle v^{\varepsilon}_N\nabla v^{\varepsilon}_N,
w^{\varepsilon}_N\rangle_{H^{0,s}}|
\end{eqnarray*}

\noindent For the first two terms in the right hand side, we apply
Lemma \ref{u,v,0,s} to obtain
$$\Big|\langle w^{\varepsilon}_N\nabla
w^{\varepsilon}_N,w^{\varepsilon}_N\rangle_{H^{0,s}}\Big|\leq
C_s\|w^{\varepsilon}_N\|_{H^{0,s}}\|\nabla_h
w^{\varepsilon}_N\|^2_{H^{0,s}},$$ and
\begin{eqnarray*}
\Big|\langle v^{\varepsilon}_N\nabla
w^{\varepsilon}_N,w^{\varepsilon}_N\rangle_{H^{0,s}}\Big|\leq
C_s\|v^{\varepsilon}_N\|^{1/2}_{H^{0,s}}\|\nabla_h
v^{\varepsilon}_N\|_{H^{0,s}}^{1/2}\|w^{\varepsilon}_N\|_{H^{0,s}}^{1/2}\|\nabla_h
w^{\varepsilon}_N\|_{H^{0,s}}^{3/2}\\
+\|w^{\varepsilon}_N\|_{H^{0,s}}\|\nabla_h
w^{\varepsilon}_N\|_{H^{0,s}}\|\nabla_h
v^{\varepsilon}_N\|_{H^{0,s}}\\
\leq
C_s\nu_h^{-3}(\nu_h^2+\|v^\varepsilon_N\|^2_{H^{0,s}})\|\nabla_h
v^\varepsilon_N\|^2_{H^{0,s}}\|
w^{\varepsilon}_N\|_{H^{0,s}}^2+\frac{\nu_h}{10}\|\nabla_h
w^{\varepsilon}_N\|_{H^{0,s}}^2.
\end{eqnarray*}

\noindent For the third term, we write
$$\langle w^{\varepsilon}_N\nabla v^{\varepsilon}_N,w^{\varepsilon}_N\rangle_{H^{0,s}}
=\langle w^{\varepsilon}_{N,h}\nabla_h
v^\varepsilon_N,w^{\varepsilon}_N\rangle_{H^{0,s}}+\langle
w^{\varepsilon}_{N,3}
\partial_3v^\varepsilon_N|w^{\varepsilon}_N\rangle_{H^{0,s}}.$$
\noindent The anisotropic product law given in Theorem 2.2 implies
\begin{eqnarray*}
\Big|\langle w^{\varepsilon}_{N,h}\nabla_h
v^\varepsilon_N,w^{\varepsilon}_N\rangle_{H^{0,s}}\Big|\leq C_s
\|\nabla_h
v^\varepsilon_N\|_{H^{0,s}}\|w^{\varepsilon}_N\|_{H^{0,s}}\|\nabla_h
w^{\varepsilon}_N\|_{H^{0,s}}\\
\leq C_s \nu_h^{-1}\|\nabla_h
v^\varepsilon_N\|^2_{H^{0,s}}\|w^{\varepsilon}_N\|^2_{H^{0,s}}+\frac{\nu_h}{10}\|\nabla_h
w^{\varepsilon}_N\|^2_{H^{0,s}}.
\end{eqnarray*}
\noindent We also have
$$
\Big|\langle
w^{\varepsilon}_{N,3}\partial_3v^\varepsilon_N,w^{\varepsilon}_N\rangle_{H^{0,s}}\Big|
\leq
C_s\|\partial_3v^\varepsilon_N\|_
 {H^{0,s}}\|w^{\varepsilon}_N\|_{H^{0,s}}\|
\nabla_h w^{\varepsilon}_N\|_{H^{0,s}}.
$$
\noindent Since $v^{\varepsilon}_N$ is localized in the $\xi_3$
direction, we can apply Lemma \ref{Bernstein} to obtain
$$\|\partial_3 v^\varepsilon_N\|_{H^{0,s}}\leq C
2^{N}\|v^\varepsilon_N\|_{H^{0,s}}.$$

\noindent Finally
\begin{eqnarray*}
\Big|\langle
w^{\varepsilon}_{N,3}\partial_3v^\varepsilon_N,w^{\varepsilon}_N\rangle_{H^{0,s}}\Big|
\leq
C_s2^N\|v^\varepsilon_N\|_{H^{0,s}}\|w^{\varepsilon}_N\|_{H^{0,s}}\|\nabla_h
w^{\varepsilon}_N\|_{H^{0,s}}\\
\leq
C_s2^{2N}\nu_h^{-1}\|v^\varepsilon_N\|^2_{H^{0,s}}\|w^{\varepsilon}_N\|^2_{H^{0,s}}+\frac{\nu_h}{10}\|\nabla_h
w^{\varepsilon}_N\|^2_{H^{0,s}}.
\end{eqnarray*}

\smallskip

\noindent What remains to estimate, therefore, is the term $\langle
v^\varepsilon_N\nabla
v^\varepsilon_N,w^{\varepsilon}_N\rangle_{H^{0,s}}$. Note that (and
this
 will be clear in the rest of the proof) this term is the only one for which we
 need that $B$ only depends on $t$ and $x_1$.\\

\noindent Since the vector field $v^\varepsilon_N$ is
divergence-free, we can write
$$
\langle v^\varepsilon_N\nabla
v^\varepsilon_N,w^{\varepsilon}_N\rangle_{H^{0,s}}=(\text{div\;}
\langle v^\varepsilon_N\otimes
v^\varepsilon_N),w^{\varepsilon}_N\rangle_{H^{0,s}}\\
= \langle v^\varepsilon_N\otimes v^\varepsilon_N,\nabla_h
w^{\varepsilon}_N\rangle_{H^{0,s}}+\langle\partial_3(v^\varepsilon_N\otimes
v^\varepsilon_N),w^{\varepsilon}_N\rangle_{H^{0,s}}.
$$

\noindent We have the estimate

\begin{eqnarray*}
\Big|\langle v^\varepsilon_N\otimes v^\varepsilon_N,\nabla_h
w^{\varepsilon}_N\rangle_{H^{0,s}}\Big| \leq
\|v^\varepsilon_N\otimes
v^\varepsilon_N\|_{H^{0,s}}\|\nabla_h w^{\varepsilon}_N\|_{H^{0,s}}\\
\leq C2^{2Ns}\|v^\varepsilon_N\otimes
v^\varepsilon_N\|_{L^2}\|\nabla_h w^{\varepsilon}_N\|_{H^{0,s}}\\
\leq
C2^{2Ns}\|v^{\varepsilon}_N\|_{L^4_{x_1}(L^\infty_{x_2,x_3})}\|v^\varepsilon_N\|_{L^4_{x_1}(L^2_{x_2,x_3})}\|\nabla_h
w^{\varepsilon}_N\|_{H^{0,
 s}}\\
\leq
C2^{(2s+1)N}\|v^\varepsilon_N\|_{L^4_{x_1}(L^2_{x_2,x_3})}^2\|\nabla_h
w^{\varepsilon}_N\|_{H^{0,s}}.
\end{eqnarray*}

\noindent By the Gagliardo-Nirenberg inequality, we get
$$
\|v^\varepsilon_N\|_{L^4_{x_1}(L^2_{x_2,x_3})}\leq
C\|v^\varepsilon_N\|_{L^2_{x_2,x_3}(L^4_{x_1})} \leq
C\|v^\varepsilon_N\|^{3/4}_{L^2}\|\partial_{x_1}v^\varepsilon_N\|_{L^2}^{1/4}.
$$

\noindent Hence

\begin{eqnarray*}
\Big|\langle v^\varepsilon_N\otimes v^\varepsilon_N, \nabla_h
w^{\varepsilon}_N\rangle_{H^{0,s}}\Big|\leq
C2^{N(2s+1)}\|v^\varepsilon_N\|_{L^2}^{3/2}\|\nabla_h
v^\varepsilon_N\|^{1/2}_{L^2}\|\nabla_h
w^{\varepsilon}_N\|_{H^{0,s}}\\
\leq C\nu_h^{-1}\|v^\varepsilon_N\|_{L^2}^3\|\nabla_h
v^{\varepsilon}_N\|_{L^2}+\frac{\nu_h}{10}\|\nabla_h
w^{\varepsilon}_N\|_{H^{0,s}}^2.
\end{eqnarray*}

\noindent Similarly, one has
\begin{eqnarray*}
\Big|\langle \partial_3(v^\varepsilon_N\otimes v^\varepsilon_N),
w^{\varepsilon}_N\rangle_{H^{0,s}}\Big|\leq
C2^{2N}2^{2Ns}\|v^\varepsilon_N
 \otimes
v^\varepsilon_N\|_{L^2}\|w^{\varepsilon}_N\|_{H^{0,s}}\\
\leq C\nu_h^{-1}\|v^\varepsilon_N\|_{L^2}^3\|\nabla_h
v^{\varepsilon}_N\|_{L^2}+\frac{\nu_h}{10}\|w^{\varepsilon}_N\|_{H^{0,s}}^2.
\end{eqnarray*}

\noindent To summarize, we combine above estimates to obtain

\begin{eqnarray*}
\frac{d}{dt}\|w^{\varepsilon}_N
(t)\|^2_{H^{0,s}}+\frac{\nu_h}{10}\|\nabla_h w^{\varepsilon}_N
(t)\|^2_{H^{0,s}} &\leq& \|w^{\varepsilon}_N
(t)\|_{H^{0,s}}\|\nabla_h w^{\varepsilon}_N
(t)\|^2_{H^{0,s}}\\&+&C\|v^\varepsilon_N(t)\|_{L^2}^3\|\nabla_h
v^\varepsilon_N(t)\|_{L^2}\\
&+&C(1+\|v^\varepsilon_N(t)\|_{H^{0,s}}^2)\|w^{\varepsilon}_N
(t)\|_{H^{0,s}}^2\\&+&C(1+\|v^\varepsilon_N(t)\|^2_{H^{0,s}})
\|\nabla_h v^\varepsilon_N(t)\|^2_{H^{0,s}}\|w^{\varepsilon}_N
(t)\|_{H^{0,s}}^2.
\end{eqnarray*}

\noindent Since $w^{\varepsilon}_N (0)$ is small in $H^{0,s}$, we
can define the time $T_{\varepsilon,N}>0$ by

$$
T_{\varepsilon,N}=\sup\Big\{\;t>0;\quad\;\forall\;0\leq t'\leq
t,\quad \|w^{\varepsilon}_N (t')\|_{H^{0,s}}\leq 2 c\nu_h\;\Big\}
$$

\noindent On the time  interval $[0,T_{\varepsilon, N}[$, we have
\begin{eqnarray*}
\frac{d}{dt}\|w^{\varepsilon}_N
(t)\|^2_{H^{0,s}}+\frac{\nu_h}{10}\|\nabla_h w^{\varepsilon}_N
(t)\|^2_{H^{0,s}} \leq C\|v^\varepsilon_N(t)\|_{L^2}^3\|\nabla_h
v^\varepsilon_N(t)\|_{L^2}\\
+C\big(1+\|v^\varepsilon_N(t)\|_{H^{0,s}}^2+(1+\|v^\varepsilon_N(t)\|^2_{H^{0,s}})\|\nabla_h
v^\varepsilon_N(t)\|^2_{H^{0,s}}\big)\|w^{\varepsilon}_N
(t)\|_{H^{0,s}}^2.
\end{eqnarray*}

\noindent Using Gronwall's Lemma, energy estimate and the fact that

$$
\dint_0^t\|\nabla_h v^\varepsilon_N(\tau)\|_{L^2}d\tau\leq
t^{1/2}\big(\dint_0^t\|\nabla_h
v^\varepsilon_N\|^2_{L^2}d\tau\big)^{1/2}\leq
\Big(\frac{t}{2\nu_h}\Big)^{1/2}\|u_0\|_{L^2},
$$
\noindent we infer
  {\small\begin{eqnarray*}
\|w^{\varepsilon}_N (t)\|^2_{H^{0,s}}\leq
\big((c\nu_h)^2+t^{\frac{1}{2}}C\|u_0\|_{L^2}^4\big)
\exp\bigg(t\;C\big(1+\|u_0\|_{H^{0,s}}^2+(1+\|u_0\|^2_{H^{0,s}})\|u_0\|^2_{L^
 2}\big)\bigg).
\end{eqnarray*}}

\noindent It follows that $T_{\varepsilon,N}\geq T$ where $T>0$ is
given by

{\small$$\big((c\nu_h)^2+T^{\frac{1}{2}}C\|u_0\|_{L^2}^4\big)
\exp\bigg(T\;C\big(\|u_0\|_{H^{0,s}}^2+(1+\|u_0\|^2_{H^{0,s}})
\|u_0\|^2_{L^2}\big)\bigg)\leq \bigg(\frac{3}{2}c\nu_h\bigg)^2.$$}

\noindent Thus we have proved that $w^{\varepsilon}_N$ exists on a
uniform time interval and the conclusion follows. This being said we
have only to add the remark that the time $T$ depends on the
distribution of the frequencies of the initial data and not only on
the size of the data.
\end{proof}

\begin{proof}[Proof of Theorem \ref{th3}]\quad\\ Let us now take
$u_0\in H^s$ with $s>1/2$. By Theorem \ref{th2} , there exists a
positive time $T$ independent of $\varepsilon>0$ and a unique
solution $u^\varepsilon\in C([0,T],H^{0,s})$ of the system
$(IRF^\varepsilon)$,  with $\nabla_h u^\varepsilon\in
L^2([0,T],H^{0,s})$. Since $H^{0,s}$ is continuously embedded in
$L^\infty_v(L^2
 _h)$, we deduce that  the norms
$\|u^\varepsilon(t)\|_{L^\infty_v(L^2_h)}$ and $\dint_0^t\|\nabla_h
u^\varepsilon(\tau)\|_{L^\infty_v(L^2_h)}^2d\tau$ are bounded on the
time interval $[0,T]$. By  Theorem \ref{thNS},  we obtain that the
lifespan $T^\varepsilon$ of the solution $u^\varepsilon$ in the
space $H^s$ satisfies the lower bound $T^\varepsilon>T$.  This means
that the solution exists on a uniform time interval in the $H^{0,s}$
space as well as in the $H^s$ space.\end{proof}


\section{\label{final} Proof of Lemma \ref{s}}

\noindent We begin with two key estimates based on the
divergence-free condition.
\begin{prop}
\label{div-free} There exists a positive constant $C$, such that for
every divergence-free vector field $u=(u_1,u_2,u_3)$,
\begin{equation}
\label{div-free-1}
 \|\nabla u_3\|_{H^s(\R^3)}\leq C\|\nabla_h
u\|_{H^s(\R^3)}
\end{equation}

\begin{equation} \label{div-fre
 e-2}
 \|\Delta_q u\|_{L^2_v
L^4_h}\leq 2^{-qs}\,c_q\,\|u\|_{H^s(\R^3)}^{1/2}\,\|\nabla_h
u\|_{H^s(\R^3)}^{1/2}.
\end{equation}

\noindent Hereafter, $(b_q)$ and $(c_q)$ stands for generic positive
sequences (which could depend on $t$) such that
$$
\dsum_{q}\,b_q\leq 1\quad{\mbox and}\quad \dsum_{q}\,c_q^2\leq 1.
$$
\end{prop}
\begin{proof}[Proof of Proposition \ref{div-free}]\quad\\
The divergence-free condition implies that
$$
\nabla u_3=\left(\nabla_h u_3, -div_h u^h\right),
$$
from which the estimate (\ref{div-free-1}) directly follows.\\

\noindent To prove (\ref{div-free-2}), we write
\begin{eqnarray*}
\|\Delta_q u\|_{L^2_v L^4_h}&\leq&\|\Delta_q
u\|_{L^2(\R^3)}^{1/2}\|\Delta_q\nabla_h u\|_{L^2(\R^3)}^{1/2}\\
&\leq&\left(2^{-\frac{qs}{2}}\,c_q^{1/2}\,\|u\|_{H^s(\R^3)}^{1/2}\right).
\left(2^{-\frac{qs}{2}}\,c_q^{1/2}\,\|\nabla_h
u\|_{H^s(\R^3)}^{1/2}\right)\\
&\leq&2^{-qs}\,c_q\,\|u\|_{H^s(\R^3)}^{1/2}\,\|\nabla_h
u\|_{H^s(\R^3)}^{1/2}.
\end{eqnarray*}
\end{proof}

\noindent Now, we return to the proof of Lemma \ref{s}. Using the
fact that
$$
\langle f,g\rangle_{H^s(\R^3)}\thickapprox\sum_{q\geq -1}\,2^{2q
s}(\Delta_q f,\Delta_q g)_{L^2(\R^3)},
$$
we only need  to prove that
\begin{eqnarray*}
\Big|(\Delta_q(u.\nabla u),\Delta_q u)_{L^2(\R^3)}\Big|&\leq& C
2^{-2qs} c_q\Big( \|u\|_{L^\infty_v L^2_h}^{1/2} \|\nabla_h
u\|_{L^\infty_v L^2_h}^{1/2}\|u\|_{H^s}^{1/2}\|\nabla_h
u\|_{H^s}^{3/2}+\\&& \|\nabla_h u\|_{L^\infty_v
L^2_h}\|u\|_{H^s}\|\nabla_h u\|_{H^s}\Big).
\end{eqnarray*}

\noindent We shall split the term $(\Delta_q(u.\nabla u),\Delta_q
u)_{L^2}$ in the following way
$$
\int_{\R^3}\Delta_q(u.\nabla u)\Delta_q u\,dx={\mathcal I}_q^h +
{\mathcal I}_q^v
$$
where
\begin{eqnarray*}
{\mathcal I}_q^h&=&\int_{\R^3}\Delta_q(u_h.\nabla_h u)\Delta_q u\, dx\\
{\mathcal I}_q^v&=&\int_{\R^3}\Delta_q(u_3\partial_3 u)\Delta_q u\,
dx.
\end{eqnarray*}

\noindent Let us start with the study of the horizontal terms. Using
the para-differential decomposition of Bony, we can write
$$
{\mathcal I}_q^h={\mathcal I}_1^{h,q}+{\mathcal I}_2^{h,q}+{\mathcal
R}^{h,q},
$$
where
\begin{eqnarray*}
{\mathcal I}_1^{h,q}&=&\sum_{|q'-q|<N_0}\,\int_{\R^3}\,
\Delta_q(S_{q'-1}u_h.\nabla_h \Delta_{q'}u)\Delta_q u\, dx\\\\
{\mathcal I}_2^{h,q}&=&\sum_{|q'-q|<N_0}\,\int_{\R^3}\,
\Delta_q(\Delta_{q'}u_h \nabla_h S_{q'-1}u )\Delta_q u\, dx\\\\
{\mathcal R}^{h,q}&=&\sum_{{\tiny{\begin{array}{c}
   q'> q-N_0\\
   i\in\{0,\pm 1\}
\end{array}}}}\int_{\R^3}\,
\Delta_q(\Delta_{q'}u_h .\nabla_h \Delta_{q'+i}u )\Delta_q u\, dx.
\end{eqnarray*}

\smallskip
\noindent$\bullet$ {\bf Estimate of ${\mathcal I}_1^{h,q}$}
\smallskip

\noindent By the H\"older inequality, we infer
$$
|{\mathcal I}_1^{h,q}|\leq\sum_{|q'-q|< N_0}
\|\Delta_q(S_{q'-1}u_h.\nabla_h \Delta_{q'}u)\|_{L^2_v
L^{4/3}_h}\|\Delta_q u\|_{L^2_v L^4_h}.
$$
Proposition \ref{div-free} gives
$$
\|\Delta_q u\|_{L^2_v L^4_h}\leq 2^{-q s}
c_q\|u\|_{H^s(\R^3)}^{1/2}\|\nabla_h u\|_{H^s(\R^3)}^{1/2}.
$$

\noindent Using Lemma \ref{uniformBound} and H\"older inequality, we
get
\begin{eqnarray*}
\|\Delta_q(S_{q'-1}u_h.\nabla_h \Delta_{q'}u)\|_{L^2_v
L^{4/3}_h}&\leq&C
\|S_{q'-1}u_h.\nabla_h \Delta_{q'}u\|_{L^2_v L^{4/3}_h}\\
&\leq&C\|S_{q'-1}u_h\|_{L^\infty_v L^4_h}\|\nabla_h
\Delta_{q'}u\|_{L^2_v L^2_h}\\
&\leq&C 2^{-q' s}c_{q'}\|u_h\|_{L^\infty_v L^4_h}\|\nabla_h u\|_{H^s(\R^3)}\\
&\leq&C 2^{-q' s}c_{q'}\|u\|_{L^\infty_v L^2_h}^{1/2}\|\nabla_h
u\|_{L^\infty_v L^2_h}^{1/2} \|\nabla_h u\|_{H^s(\R^3)}.
\end{eqnarray*}

\noindent Finally
$$
|{\mathcal I}_1^{h,q}|\leq 2^{-2q s}b_q\|u\|_{L^\infty_v
L^2_h}^{1/2} \|\nabla_h u\|_{L^\infty_v
L^2_h}^{1/2}\|u\|_{H^s(\R^3)}^{1/2} \|\nabla_h
u\|_{H^s(\R^3)}^{3/2}.
$$

\smallskip

\noindent$\bullet$ {\bf Estimate of ${\mathcal I}_2^{h,q}$}

\smallskip

\noindent Arguing as in the estimate of the term ${\mathcal
I}_1^{h,q}$, we  get
\begin{eqnarray*}
|{\mathcal I}_2^{h,q}|&\leq&\sum_{|q- q'|< N_0}
\|\Delta_q(\Delta_{q'}u_h.\nabla_h S_{q'-1}
 u)\|_{L^2_v
L^{4/3}_h}\|\Delta_q u\|_{L^2_v L^4_h}\\
&\leq&C\sum_{|q-q'|< N_0}\|\Delta_{q'} u\|_{L^2_v L^4_h}\|\nabla_h
S_{q'-1}u\|_{L^\infty_v L^2_h}
\|\Delta_q u\|_{L^2_v L^4_h}\\
&\leq&2^{-2 q s}b_q\|\nabla_h u\|_{L^\infty_v
L^2_h}\|u\|_{H^s}\|\nabla_h u\|_{H^s}.
\end{eqnarray*}

\smallskip

\noindent$\bullet$ {\bf Estimate of ${\mathcal R}^{h,q}$}

\smallskip

\noindent Using the H\"older inequality, we get
$$
|{\mathcal R}^{h,q}|\leq \|\Delta_q u\|_{L^2_v
L^4_h}\sum_{{\tiny{\begin{array}{c}
   q'\geq q-4\\
   i\in\{0,\pm 1\}
\end{array}}}}\|\Delta_q(\Delta_{q'}u_h .\nabla_h \Delta_{q'+i}u )\|_{L^2_v
L^{4/3}_h}
$$
\noindent Proposition \ref{div-free} implies that
$$
\|\Delta_q u\|_{L^2_v L^4_h}\leq 2^{-q s}
c_q\|u\|_{H^s(\R^3)}^{1/2}\|\nabla_h u\|_{H^s(\R^3)}^{1/2}.
$$
\noindent Moreover, we have
\begin{eqnarray*}
\|\Delta_q(\Delta_{q'}u_h .\nabla_h \Delta_{q'+i}u )\|_{L^2_v
L^{4/3}_h}&\leq&C\|\Delta_{q'}u_h .\nabla_h \Delta_{q'+i}u \|_{L^2_v
L^{4/3}_h}\\&\leq&C\|\Delta_{q'}u_h \|_{L^\infty_v
L^{4}_h}\|\nabla_h \Delta_{q'+i}u\|_{L^2(\R^3)}\\&\leq&2^{-
q's}c_{q'}\|\nabla_h u\|_{H^s(\R^3)}\|u\|_{L^\infty_v
L^{2}_h}^{1/2}\|\nabla_h u\|_{L^\infty_v L^{2}_h}^{1/2}.
\end{eqnarray*}

\noindent Finally, we obtain
$$
|{\mathcal R}^{h,q}|\leq 2^{-2qs}b_q\|u\|_{L^\infty_v L^2_h}^{1/2}
\|\nabla_h u\|_{L^\infty_v L^2_h}^{1/2}\|u\|_{H^s}^{1/2}\|\nabla_h
u\|_{H^s}^{3/2}.
$$

\smallskip

\noindent$\bullet$ {\bf Estimate of ${\mathcal I}^v _q$}

\smallskip

\noindent The vertical term is more delicate since the vertical
viscosity vanishes, and hence we loose the regularity in the
vertical direction. To avoid this difficulty, we proceed in the
following way. We write  the Bony's decomposition for $\Delta_q(u_3\partial_3 u)$ in two para-product terms and in the rest term. Using some commutators, we decompose now the para-product term $\Delta_qT_{u_3}\partial_3u$ as following
\begin{equation*}
\Delta_q(S_{q'-1}u_3\partial_3\Delta_{q'}u)=S_{q-1}u_3\partial_3
 \Delta_q \Delta_{q'}u+(S_{q'-1}u_3-S_{q-1}u_3)\partial_3\Delta_q\Delta_{q'}u+[\Delta_q, S_{q'-1}u_3]\partial_3\Delta_{q'}u.
\end{equation*}
Summing in $q'$, we obtain finally the following decomposition
$$
{\mathcal I}^v _q={\mathcal I}_1^{v,q} +{\mathcal
I}_2^{v,q}+{\mathcal I}_3^{v,q}+{\mathcal I}_4^{v,q}+ {\mathcal
R}^{v,q},
$$
where
\begin{eqnarray*}
{\mathcal I}_1^{v,q}&=&\int_{\R^3}\,(S_{q-1}u_3)\partial_3\Delta_q
u.\Delta_q u\,dx \\\\
{\mathcal I}_2^{v,q}&=&\sum_{|q'-q|<
N_0}\int_{\R^3}\,[\Delta_q,S_{q'-1}u_3]
\partial_3\Delta_{q'}u.\Delta_q u\,dx \\\\
{\mathcal I}_3^{v,q}&=&\sum_{|q'-q|< N_0}\int_{\R^3}\,(S_{q'-1}
u_3-S_{q-1}u_3)
\partial_3\Delta_q\Delta_{q'}u.\Delta_q u\,dx \\\\
{\mathcal I}_4^{v,q}&=&\sum_{|q'-q|<
N_0}\int_{\R^3}\,\Delta_q\left(\Delta_{q'}u_3
\partial_3 S_{q'-1}u\right).\Delta_q u\,dx \\\\
{\mathcal R}^{v,q}&=&\sum_{{\tiny{\begin{array}{c}  q'>q-N_0\\\\
|q''-q'|\leq 1\end{array}}}}
\int_{\R^3}\,\Delta_q\left(\Delta_{q'}u_3\partial_3\Delta_{q''}u\right).
\Delta_q
u\,dx.
\end{eqnarray*}

\smallskip

\noindent We shall now estimate each of the above terms.

\smallskip

\noindent$\bullet${\bf Estimate of ${\mathcal I}_1^{v,q}$}

\smallskip

\noindent Integrating by parts yields
$$
{\mathcal I}_1^{v,q}=-\frac{1}{2}\int_{\R^3}\,(S_{q-1}\partial_3
u_3)\Delta_q u. \Delta_q u\,dx.
$$
\noindent Using Proposition \ref{div-free}, we get
\begin{eqnarray*}
|{\mathcal I}_1^{v,q}|&\leq&C\|\nabla_h u\|_{L^\infty_v
L^2_h}\|\Delta_q u\|_{L^2_v L^4_h}^2\\
&\leq& 2^{-2 qs} b_q\|\nabla_h u\|_{L^\infty_v
L^2_h}\|u\|_{H^s}\|\nabla_h u\|_{H^s}.
\end{eqnarray*}

\smallskip

\noindent$\bullet${\bf Estimate of ${\mathcal I}_2^{v,q}$}

\smallskip

\noindent By H\"older inequality, we have
$$
|{\mathcal I}_2^{v,q}|\leq\|[\Delta_q,S_{q'-1}u_3]
\partial_3\Delta_{q'}u\|_{L^2_v L^{4/3}_h}\|\Delta_q u\|_{L^2_v L^4_h}.
$$
Using Lemma \ref{commutator}, Proposition \ref{div-free} and the
fact that $\Delta_{q'}u$ is localized in vertical frequencies, we
infer
\begin{eqnarray*}
|{\mathcal I}_2^{v,q}|&\leq& C 2^{-q}\|S_{q'-1}\nabla
u_3\|_{L^\infty_v L^2_h}\|\partial_3
\Delta_{q'}u\|_{L^2_v L^4_h}\|\Delta_q u\|_{L^2_v L^4_h}\\
&\leq&C 2^{q'-q}\,\,\|(\nabla_h u_3,\partial_3 u_3)\|_{L^\infty_v
L^2_h}\,\|
\Delta_{q'}u\|_{L^2_v L^4_h}\|\Delta_q u\|_{L^2_v L^4_h}\\
&\leq&C\|\nabla_h u\|_{L^\infty_v L^2_h}\|
\Delta_{q'}u\|_{L^2_v L^4_h}\|\Delta_q u\|_{L^2_v L^4_h}\\
&\leq&C\|\nabla_h u\|_{L^\infty_v L^2_h}\|\Delta_{q'}u\|_{L^2}^{1/2}
\|\Delta_{q'}\nabla_h u\|_{L^2}^{1/2}\|\Delta_{q}u\|_{L^2}^{1/2}
\|\Delta_{q}\nabla_h u\|_{L^2}^{1/2}\\
&\leq&C2^{-2qs}b_q\|\nabla_h u\|_{L^\infty_v
L^2_h}\|u\|_{H^s}\|\nabla_h u\|_{H^s}.
\end{eqnarray*}

\smallskip

\noindent$\bullet${\bf Estimate of ${\mathcal I}_3^{v,q}$}

\smallskip

\noindent Since $S_{q'-1}u_3-S_{q-1}u_3$ is localized in frequency
in a ring
  of size $2^q$, we get
\begin{eqnarray*}
|{\mathcal I}_3^{v,q}|&\leq&\sum_{q'\sim
q}\|S_{q'-1}u_3-S_{q-1}u_3\|_{L^\infty_v L^2_h}
\|\partial_3
 \Delta_q\Delta_{q'}u\|_{L^2_v L^4_h}\|\Delta_q u\|_{L^2_v L^4_h}\\
&\leq&\sum_{q'\sim q}2^{-q}\|\big(S_{q'-1}-S_{q-1}\big)\nabla
u_3\|_{L^\infty_v L^2_h} \|\partial_3\Delta_q\Delta_{q'}u\|_{L^2_v
L^4_h}\|\Delta_q
u\|_{L^2_v L^4_h}\\
\end{eqnarray*}

\noindent Using Proposition \ref{div-free}, we obtain
\begin{eqnarray*}
 |{\mathcal I}_3^{v,q}|&\leq&C
2^{-2qs}b_q\|\nabla_h u\|_{L^\infty_v L^2_h}\|u\|_{H^s}\|\nabla_h
u\|_{H^s}.
\end{eqnarray*}

\smallskip

\noindent$\bullet${\bf Estimate of ${\mathcal I}_4^{v,q}$}

\smallskip

\noindent As in the above estimates, we have
\begin{eqnarray*}
|{\mathcal I}_4^{v,q}|&\leq& C\sum\limits_{|q-q'|<
N_0}\|\Delta_q\left(\Delta_{q'}u_3\partial_3S_{q'-1}u\right)\|_
{L^2_v L^{4/3}_h}\|\Delta_q u\|_{L^2_v L^4_h}\\
&\leq&C\sum\limits_{|q-q'|<
N_0}\|\Delta_{q'}u_3\partial_3S_{q'-1}u\|_
{L^2_v L^{4/3}_h}\|\Delta_q u\|_{L^2_v L^4_h}\\
&\leq&C\sum\limits_{|q-q'|< N_0}\|\Delta_{q'}u_3\|_{L^2_v
L^2_h}\|\partial_3 S_{q'-1}u\|_{L^\infty_v L^4_h} \|\Delta_q
u\|_{L^2_v L^4_h}.
\end{eqnarray*}

\noindent Proposition \ref{div-free} implies that
\begin{eqnarray*}
\|\Delta_{q'}u_3\|_{L^2}&\leq&C 2^{-q'}\|\Delta_{q'}\nabla u_3\|_{L^2}\\
&\leq&C 2^{-q}2^{-q' s} c_{q'}\|\nabla_h u\|_{H^s},
\end{eqnarray*}
and
$$
\|\Delta_q u\|_{L^2_v L^4_h}\leq
2^{-qs}c_q\|u\|_{H^s}^{1/2}\|\nabla_h u\|_{H^s}^{1/2}.
$$

\noindent In addition
\begin{eqnarray*}
\|\partial_3 S_{q'-1}u\|_{L^\infty_v L^4_h}&\leq&C
2^{q'}\|u\|_{L^\infty_v L^4_h}\\
&\leq&C 2^{q'}\|u\|_{L^\infty_v L^2_h}^{1/2}\|\nabla_h
u\|_{L^\infty_v L^2_h}^{1/2}.
\end{eqnarray*}

Hence
$$
|{\mathcal I}_4^{v,q}|\leq 2^{-2q s}b_q\|u\|_{L^\infty_v
L^2_h}^{1/2} \|\nabla_h u\|_{L^\infty_v
L^2_h}^{1/2}\|u\|_{H^s}^{1/2}\|\nabla_h u\|_{H^s}^{3/2}.
$$

\smallskip

\noindent$\bullet${\bf Estimate of ${\mathcal R}^{v,q}$}

\smallskip

\noindent We easily see that
\begin{eqnarray*}
|{\mathcal
R}^{v,q}|&\leq&C\sum\limits_{q'>q-N_0}\|\Delta_q\left(\Delta_{q'}u_3\partial_3\Delta_{q''}u\right)\|_{L^
 2_v
L^{4/3}_h}\|\Delta_q u\|_{L^2_v L^4_h}\\
&\leq&C\sum\limits_{q'>q-N_0}\|\Delta_{q'}u_3\partial_3\Delta_{q''}u\|_{L^2_v
L^{4/3}_h}\|\Delta_q u\|_{L^2_v L^4_h}\\
&\leq&\sum\limits_{q'>q-N_0}\|\Delta_{q'}u_3\|_{L^2}\|\partial_3\Delta_{q''}u\|_{L^\infty_v
L^4_h} \|\Delta_q u\|_{L^2_v L^4_h}.
\end{eqnarray*}

\noindent Using Proposition \ref{div-free} and Lemma
\ref{Bernstein}, we infer
\begin{eqnarray*}
\|\Delta_{q'}u_3\|_{L^2}&\leq&2^{-q'}\|\Delta_{q'}\nabla u_3\|_{L^2}\\
&\leq&2^{-q'}\|\Delta_{q'}\nabla_h u_3\|_{L^2}\\
&\leq&2^{-q'}2^{-q' s}c_{q'}\|\nabla_h u\|_{H^s},
\end{eqnarray*}
\begin{eqnarray*}
\|\partial_3\Delta_{q''}u\|_{L^\infty_v L^4_h}&\leq&
2^{q''}\|u\|_{L^\infty_v L^4_h}\\
&\leq&2^{q''}\|u\|_{L^\infty_v L^2_h}^{1/2}\|\nabla_h
u\|_{L^\infty_v L^2_h}^{1/2},
\end{eqnarray*}
and
$$
\|\Delta_q u\|_{L^2_v L^4_h}\leq 2^{-qs}
c_q\|u\|_{H^s}^{1/2}\|\nabla_h u\|_{H^s}^{1/2}.
$$

\noindent Finally
$$
|{\mathcal R}^{v,q}|\leq 2^{-2qs}b_q\|u\|_{L^\infty_v L^2_h}^{1/2}
\|\
 nabla_h u\|_{L^\infty_v L^2_h}^{1/2}\|u\|_{H^s}^{1/2}\|\nabla_h
u\|_{H^s}^{3/2}.
$$
The proof of Lemma \ref{s} is completed. \endproof


\section{Appendix}

\noindent This section is devoted to the proofs of some technical
lemmas stated in the second section. Note that the proofs of such
results are contained in several papers (see for instance
\cite{chem, paicu, paicu1}). Here we give them for the convenience
of the reader.

\smallskip

\begin{proof}[Proof of Bernstein Lemma \ref{anablaa}]\quad\\
Let $\tilde{\mathcal C}$ be a ring such that ${\mathcal
C}\subset\tilde{\mathcal C}$ and $\phi$ be a smooth compactly
supported function in $\tilde{\mathcal C}$, whose value is $1$ near
${\mathcal C}$. Since $\text{supp}\,{\mathcal F}\,u\subset 2^q
{\mathcal C}$, we infer
$$
{\mathcal F}\, \triangle_q u(\xi)=\phi(2^{-q}\xi){\mathcal F}\,
\triangle_q
u(\xi)=\sum\limits_{i=1,3}\frac{\phi_i(2^
 {-q}\xi)}{\xi_i}\widehat{
\partial_{x_i}\triangle_q u}(\xi),
$$
where $\phi(\xi)=\sum\limits_{i=1,3}\phi_i(\xi)$ is such that
$\xi_i\neq 0$ on the support of $\phi_i$. If we set
$\psi_i(\xi)=~\phi_i(\xi)/{\xi_i}$, then $\psi_i\in
C^\infty_0(\R^3)$, has the support in $\mathcal C$ and

$${\mathcal
F}\,\triangle_qu(\xi)=2^{-q}\sum\limits_{i=1,3}\psi_i(2^{-q}\xi)\widehat{
\partial_{x_i}\triangle_q u}(\xi).$$
We denote by $h$ the rapidly decaying  function such that ${\mathcal
F}\,h_i =\psi_i$.  Therefore
$$
\triangle_qu(x)=2^{-q}2^{3q}\sum\limits_{i=1,3}h_i(2^q.)\ast\partial_{x_i}
\triangle_q u(x).$$

\noindent From this equality, we easily deduce that
$$
\|\triangle_qu\|_{L^p_h(L^r_v)}\leq C
2^{-q}\sum\limits_{i=1,3}\|2^{3q}h_i(2^q\cdot)\|_{L^1}\|\partial_{x_i}\triangle_q
u\|_{L^p_h(L^r_v)}
$$
and so
$$
\|\triangle_q u\|_{L^r_v(L^p_h)}\leq C 2^{-q}\|\nabla \triangle_q
u\|_{L^r_v(L^p_h)}.
$$
\end{proof}

\

\begin{proof}[Proof of Lemma \ref{uniformBound}]\quad\\
\noindent We introduce the function $h_j(x)=2^{3j}h(2^jx)$, where
$h\in\mathcal S(\R^3)$ is such that $\mathcal F (h)=\varphi\in
C_0^\infty(\R^3)$. Let
  \begin{eqnarray*}
  u_j(x):=\Delta_j u(x)&=&\int_{\R^3}\,h_j(x-y) u(y) dy\\
  &=&\int_{\R^3}\,h_j(x_h-y_h,x_v-y_v) u(y_h,y_v) dy_h d y_v\\
  &=&\int_{\R}\left(\int_{\R^2}\,h_j(x_h-y_h,x_v-y_v) u(y_h,y_v) dy_h\right) d y_v\\
  &=&\int_{\R}\big(h_j(\cdot,x_v-y_v)\ast_{x_h} u(\cdot,y_v)\big)(x_h) d
  y_v.
  \end{eqnarray*}

  \noindent Then
  $$
  \|u_j(\cdot,x_v)\|_{L^r_h}\leq\int_{\R}\,\left\|h_j(\cdot, x_v-y_v)\ast_{x_h}
u(\cdot,y_v)\right\|_{L^r_h}dy_v.
  $$

  \noindent By the Young inequality, we get
  $$
\left\|h_j(\cdot,x_v-y_v)\ast_{x_h} u(\cdot,y_v)\right\|_{L^r_h}\leq
\|h_j(\cdot,x_v-y_v)\|_{L^1_h}
  \|u(\cdot, y_v)\|_{L^r_h}
  $$
  and so
  \begin{eqnarray*}
\|u_j(\cdot,x_v)\|_{L^r_h}&\leq&\int_{\R}\|h_j(\cdot,x_v-y_v)\|_{L^1_h}
  \|u(\cdot, y_v)\|_{L^r_h}dy_v\\
  &=&\left(\|h_j(.,..)\|_{L^1_h}\ast_{x_v}\|u(.,..)\|_{L^r_h}\right)
 (x_v).
  \end{eqnarray*}

  \noindent Now, we take the $L^p_v$ norm in the above inequality to obtain
  $$
  \Big\|\|u_j(\cdot,x_v)\|_{L^r_h}\Big\|_{L^p_v}\leq\Big\|\|h_j(.,..)\|_{L^1_h}\Big\|_{L^1_v}
\Big\|\|u(.,..)\|_{L^r_h}\Big\|_{L^p_v}.
$$

\noindent Finally
$$
\|\Delta_j u\|_{L^p_v L^r_h}\leq \|h\|_{L^1(\R^3)}\|u\|_{L^p_v
L^r_h}.
$$
\end{proof}

\begin{proof}[Proof of Lemma \ref{commutator}]\quad\\
\noindent We will prove here, that, for all $p,r,s,t\geq 1$ such
that $\frac{1}{r}=\frac{1}{s}+\frac{1}{t}$, we have

$$\|[\Delta_q; a]b\|_{L^p_v L^r_h}\leq C 2^{-q}\|\nabla
a\|_{L^\infty_v L^s_h}\|b\|_{L^p_v L^t_h}.$$

\noindent We begin by writing the commutator in the following form

\begin{eqnarray*}
[\Delta_q;a]b(x)&=&2^{3q}\int\limits_{\R^3}h(2^q y)(a(x-y)-a(x))b(x-y)dy\\
&=&-2^{3q}\int\limits_{\R^3}h(2^q y)\int_0^1y\nabla a(x-\tau y)d\tau b(x-y)dy\\
&=&-2^{3q}\int\limits _{\R^3\times[0,1]}\sum\limits_{i=1,3}h(2^q
y)y_i(\partial_i a)(x-\tau y)b(x-y)dyd\tau\\
&=&2^{
 -q}2^{3q}\int\limits_{\R^3\times[0,1]}\sum\limits_{i=1,3}h_i(2^q
y)\cdot (\partial_i a)(x- \tau y)b(x-y)dyd\tau,
\end{eqnarray*}
\noindent where $h_i(z)=h(z)z_i$. We take the $L^r_h$ norm in the
horizontal variable and we apply the H\"older inequality to obtain
{\small
\begin{eqnarray*}
\|[\Delta_q;a]b(\cdot,x_v)\|_{L^r_h}&\leq&2^{-q}
\int\limits_{\R^3\times[0,1]}\sum\limits_{i=1,3}2^{3q}|h_i(2^q
y)|\|(\partial_i
a)(\cdot-\tau y_h,x_v-\tau y_v)b(\cdot-y_h,x_v-y_v)\|_{L^r_{x_h}}dyd\tau\\
&\leq& 2^{-q}\int\limits_{\R^3\times
[0,1]}\sum\limits_{i=1,3}2^{3q}|h_i(2^q y)| \|(\partial_i a)(\cdot-
\tau y_h,x_v-\tau y_v )\|_{L^s_
{x_h}}\|b(\cdot-y_h,x_v-y_v)\|_{L^t_{x_h}}dyd\tau\\
&=&2^{-q}\int\limits_{\R^3\times
[0,1]}\sum\limits_{i=1,3}2^{3q}|h_i(2^q y)| \|\partial_i a(\cdot,
x_v-\tau y_v)\|_{L^s_ h}\|b(\cdot,x_v-y_v)\|_{L^t_h}dyd\tau.
\end{eqnarray*}
}

\noindent We take now the $L^p_v$ norm in the vertical variable and
we obtain

{\small
\begin{eqnarray*}
\|[\Delta_q;a]b\|_{L^p_v
 L^r_h} &\leq&
2^{-q}\int\limits_{\R^3\times[0,1]}\sum\limits_{i=1,3}2^{3q}|h_i(2^q
y)| \big\|\|(\partial_i a)(\cdot,x_v-\tau
y_v)\|_{L^s_h}\|b(\cdot, x_v-y_v)\|_{L^t_h}\big\|_{L^p_{x_v}}dyd\tau\\
&\leq& 2^{-q}\int\limits_{\R^3\times
[0,1]}\sum\limits_{i=1,3}2^{3q}|h_i(2^q y)|
\sup\limits_{x_v}\|\partial_i a(\cdot,x_v-\tau y_v)\|_
{L^s_{x_h}}\|\|b(\cdot, x_v-y_v)\|_{L^t_{x_h}}\|_{L^p_{x_v}}dyd\tau\\
  &=&2^{-q}\int\limits_{\R^3\times
[0,1]}\sum\limits_{i=1,3}2^{3q}|h_i(2^q y)| \|\partial_i
a\|_{L^\infty_v(L^s_h)}
  \|b\|_{L^p_v(L^t_h)}dyd\tau\\
  &=&2^{-q}\|\nabla a\|_{L^\infty_v
L^s_h}\|b\|_{L^p_v(L^t_h)}\|h_i\|_{L^1(\R^3)}.
\end{eqnarray*}
}

\noindent This leads to
$$
\|[\Delta_q; a]b\|_{L^p_v(L^r_h)}\leq C 2^{-q}\|\nabla
a\|_{L^\infty_v(L^s_h)}\|b\|_{L^p_vL^t_h}.
$$
\end{proof}
\noindent{\bf Acknowledgement}\\

\noindent {\sf We are deeply indebted to the referee for his
suggestions
  which helped us to improve this paper}.



\begin{thebibliography}{10}






\bibitem{BMN}
{\bf A.Babin, A. Mahalov and B.Nicolaenko}, \newblock{\em Regularity
and integrability of 3D Euler and Navier-Stokes equations for
rotating fluids}, Asymptotic Analysis, 15, p. 103-150, 1997.






\bibitem {BMN1}
{\bf A.Babin, A. Mahalov and B.Nicolaenko}, \newblock{\em Global
Regularity of 3D Rotating Navier-Stokes Equations for Resonant
Domains}, Indiana University Mathematics Journal, Vol. 48, No. 3, p.
1133-1176, 1999.



\bibitem {bony}
{\bf J.-M. Bony}, \newblock{\em Calcul symbolique et propagation des
singularit\'es pour les \'equations aux d\'eriv\'ees partielles non
lin\'eaires}, Annales scientifiques de l'\'Ecole Normale
Sup\'erieure, 14, p. 209-246, 1981.



\bibitem {bresch}
{\bf D. Bresch, D. G\'erard-Varet and E. Grenier},
\newblock{\em Derivation of the Planetary Geostrophic E
 quations},
Arch. Rat. Mech. Anal. {\bf 182}, No 2, p. 387-413, 2006.



\bibitem {chem}
{\bf J.-Y. Chemin},  \newblock{\em Localization in Fourier space and
Navier-Stokes system}, Preprint.


\bibitem {chem1}
{\bf J.-Y. Chemin},  \newblock{\em Fluides parfaits incompessibles},
Ast\'erisque {\bf 230}, 1995.








\bibitem{CDGG}
{\bf J.-Y. Chemin, B. Desjardins, I. Gallagher and E. Grenier},
\newblock{\em Fluids with anisotropic
viscosity},  M2AN. Math. Numer. Anal., 34, No. 2, p. 315-335, 2000.

\bibitem{CDGG1}
{\bf J.-Y. Chemin, B. Desjardins, I. Gallagher and E. Grenier},
\newblock{\em Ekman boundary layers in rotating fluids}, ESAIM Contr\^ole
optimal et Calcul des Variatuions, Special Tribute issue to
Jacques-Louis Lions, {\bf 8}, p. 441-466, 2002.



\bibitem{Cheverry}
{\bf C. Cheverry}, \newblock{\em Propagation of oscillations in Real
Vanishing Viscosity Limit}, Commun. Math. Phys. 247, p. 655-695,
2004.






\bibitem{FK}
{\bf H. Fujita and T. Kato}, \newblock{\em O
 n the Navier-Stokes
initial value problem I}, Archiv for Rational Mechanic and Analysis,
16, p. 269-315, 1964.

\bibitem{G}
{\bf I. Gallagher}, {\it Applications of Schochet's Methods to
Parabolic Equation}, Journal de Math\'ematiques Pures et
Appliqu\'ees, 77, p. 989-1054, 1998.

\bibitem{G1}
{\bf I. Gallagher},
  {\it Asymptotics of the Solutions of Hyperbolic
Equations With a Skew-Symmetric Perturbation}, Journal of
Differential Equations, 150, p. 363-384, 1998.










\bibitem{GR}
{\bf I. Gallagher and L. Saint-Raymond},
  \newblock{\em Weak convergence results for inhomogeneous rotating
fluid equations}, Journal d'Analyse Math\'ematique, Vol. 99, p.
1-34, 2006.

\bibitem{GSR}
{\bf I. Gallagher and L. Saint-Raymond},
  \newblock{\em On the influence of the Earth's rotation on geophysical flows},
  Handbook of Mathematical Fluid Dynamics Vol 4, p. 201-329.

\bibitem{GSR1}
{\bf I. Gallagher and Laure Saint-Raymond},
  \newblock{\em On pressureless gases driven by a stro
 ng inhomogeneous magnetic field},
  SIAM Journal for Mathematical Analysis,  36,  no. 4, p. 1159-1176, 2005.

   \bibitem{Gr}
{\bf E. Grenier},
  {\it Oscillatory Perturbations of the Navier-Stokes
Equations}, Journal de Math\'ematiques Pures et Appliqu\'ees, 76, p.
477-498, 1997.

\bibitem{GM}
{\bf E. Grenier and N. Masmoudi},  \newblock{\em Ekman layers of
rotating fluids, the case of well
 prepared initial data},  Comm. Partial Differential Equations  22,
 no. 5-6, p. 953-975, 1997.


\bibitem{Iftimie}
{\bf D. Iftimie}, \newblock{\em Resolution of the  Navier-Stokes
equations in anisotropic spaces}, Revista Mat\'ematica
Iberoamericana, 15, no. 1, p. 1-36, 1999.





\bibitem{Iftimie1}
{\bf D. Iftimie}, \newblock{\em The 3D Navier-Stokes equations seen
as a perturbation of the 2D Navier-Stokes equations}, Bulletin de la
Soci\'et\'e Math\'ematique de France, 127, p. 473-517, 1999.






\bibitem{Iftimie2}
{\bf D. Iftimie}, \newblock{\em A Uniqueness result for the
Navier-St
 okes Equations with
   vanishing vertical viscosity},  SIAM J. Math. Anal. Vol. 33, No. 6,
p. 1483-1493, 2002.


\bibitem{JMR1} {\bf J. Joly, G. M\'etivier and J. Rauch},
\newblock{\em Generic rigourous Asymptotic Expansions for Weakly Nonlinear
Multidimensional Oscillatory Waves}, Duke Mathematical Journal, {\bf
70}, p. 373-404, 1993.

\bibitem{JMR2} {\bf J. Joly, G. M\'etivier and J. Rauch},
\newblock{\em Coherent and Focusing Multidimensional Oscillatory Nonlinear
Geometric optics}, Annales Scientifiques de L'ENS, {\bf 28}, p.
51-113, 1995.


\bibitem{KM} {\bf S. Klainerman and A. Majda},
\newblock{\em Singular Limits of Quasilinear Hyperbolic system with Large
Parameters, and the Incompressible Limit of Compressible Fluids},
Communications on pure and applied Mathematics, {\bf 34}, p.
481-524, 1981.









\bibitem{leray}
{\bf J. Leray}, \newblock{\em Sur le mouvement d'un liquide visqueux
remplissant l'espace}, Acta Math, 63, p. 193-248, 1934.



\bibitem{Mas1} {\bf N
 . Masmoudi},
\newblock{\em Ekman layers of rotating fluids: the case of
general initial data}, Comm. Pure Appl. Math.  53,  no. 4, p.
432-483, 2000.





\bibitem{paicu}
{\bf M. Paicu},  \newblock{\em \'Equation  anisotrope de
Navier-Stokes dans des espaces critiques },  Rev. Mat.
Iberoamericana, 21, no. 1, p. 179--235, 2005.

\bibitem{paicu1}
{\bf M. Paicu}, \newblock{\em \'Etude asymptotique pour les fluides
anisotropes en rotation rapide dans le cas p\'eriodique}, Journal
des Math\'ematiques Pures et Appliqu\'ees, 83, p. 163-242, 2004.

\bibitem{paicu2}
{\bf M. Paicu}, \newblock{\em \'Equation p\'eriodique de
Navier-Stokes sans viscosit\'e dans une direction},  Comm. Partial
Differential Equations,  30, no. 7-9, p. 1107--1140, 2005.





\bibitem{ST}
{\bf M. Sabl\'e-Tourgeron}, \newblock{\em R\'egularit\'e microlocale
pour des probl\`emes aux limites non lin\'eaires}, Ann. Inst.
Fourier, 36, p. 39-82, 1986.

\bibitem{S1} {\bf S. Schochet},
\newblock{\em The Compressible E
 uler Equations in a Bounded Domain: Existence
of Solutions and the Incompressible Limit}, Comm. Math. Physics,
{\bf 104}, p. 49-75, 1986.

\bibitem{S2} {\bf S. Schochet},
\newblock{\em Fast Singular Limits of Hyperbolic PDEs}, Journal of
Differential Equations, {\bf 114}, p. 476-512, 1994.















\end{thebibliography}
\end{document}